\documentclass[twoside,reqno]{amsart}

\usepackage{amsmath}
\usepackage{amssymb,latexsym,amstext,amsthm}
\usepackage{verbatim} 
\usepackage{eucal} 
\usepackage{color}
\theoremstyle{plain}
\newtheorem{thm}[equation]{Theorem}
\newtheorem{pro}[equation]{Proposition}
\newtheorem{cor}[equation]{Corollary}
\newtheorem{lem}[equation]{Lemma}

\theoremstyle{definition}

\newtheorem{DEF}[equation]{Definition}
\newtheorem{rem}[equation]{Remark}

\def\rdot{\dot R}
\def\rdotcross{\rdot^\times}
\def\rcross{R^\times}
\def\op{\oplus}
\def\andd{\quad\hbox{and}\quad}

\def\RR{\mathcal R}

\def\PP{\mathcal P}

\def\MM{\mathcal M}

\def\Ab{\bar A}
\def\a{\alpha}
\def\ab{\bar\a}
\def\bb{\bar\b}

\def\w{{\mathcal W}}
\def\rd{\dot{R}}

\def\la{\langle}
\def\ra{\rangle}

\def\d{\delta}
\def\b{\beta}
\def\db{\dot\b}

\def\sg{\sigma}
\def\bbbz{{\mathbb Z}}

\def\da{\dot{\a}}
\def\dg{\dot{\gamma}}
\def\rank{\hbox{rank}}

\def\Aut{\hbox{Aut}}
\def\bd{\dot{\b}}

\def\ep{\epsilon}
\def\fm{(\cdot,\cdot)}
\def\sub{\subseteq}
\def\rds{\dot{R}_{sh}}
\def\rdl{\dot{R}_{lg}}

\def\u{{\mathcal U}}
\def\v{{\mathcal V}}
\def\w{{\mathcal W}}
\def\r{R}
\def\vb{\overline{\v}}
\def\supp{\hbox{supp}}

\def\rds{\dot{R}_{sh}}
\def\rdl{\dot{R}_{lg}}
\def\rde{\dot{R}_{ex}}
\def\bbbq{{\mathbb Q}}

\def\rlg{R_{lg}}
\def\rsh{R_{sh}}
\def\dgamma{\dot{\gamma}}
\def\ot{\otimes}
\def\Rb{\bar{R}}

\begin{document}

\setcounter{page}{1} \setcounter{page}{1}

\author{Saeid Azam$^{1,\dag}$, Hiroyuki Yamane$^{2}$, Malihe Yousofzadeh$^{1,\dag}$}

\title{Reflectable bases for affine reflection systems}

\dedicatory{Dedicated to the memory of Valiollah Shahsanaei}





\thanks{$\;^1$Department of Mathematics, University of Isfahan, Isfahan, Iran,
P.O.Box 81745-163 and School of Mathematics, Institute for
Research in Fundamental Sciences (IPM), P.O. Box: 19395-5746,
Tehran, Iran, azam@sci.ui.ac.ir, ma.yousofzadeh@sci.ui.ac.ir}

\thanks{$\;^2$Department of Pure and Applied Mathematics, Graduate
School of Information Science and Technology, Osaka University,
Toyonaka, Osaka, 560-0043, Japan, yamane@math.sci.osaka-u.ac.jp.}

\thanks{$\;^\dag$This research was in part supported by a grant from IPM
(No. 89170216, No. 89170030). The authors would   like to thank
the Center of Excellence for Mathematics, University of Isfahan}

\begin{abstract} The notion of a ``root base'' together with its
geometry plays a crucial role in the theory of finite and affine
Lie theory. However, it is known that such a notion does not exist
for the recent generalizations of finite and affine root systems
such as extended affine root systems and affine reflection
systems. In this work, we consider the notion of a ``reflectable base`` for an affine reflection system $R.$ A reflectable base for $R$ is  a minimal subset $\Pi$ of roots such that
the non-isotropic part of the root system can be recovered by
reflecting roots of $\Pi$ relative to the hyperplanes determined
by $\Pi$. We give a full characterization of reflectable bases for
tame irreducible affine reflection systems of reduced types,
excluding types $E_{6,7,8}$. As a by-product of our results, we
show that if the root system under consideration is locally finite,
then any reflectable base is an integral base.
\end{abstract}

%
%
%

\maketitle

\parbox{4.5in}{\small\it2000 Mathematics Subject
Classification. 17B67, 17B65, 20F55, 51F15.

 }

\medskip
\setcounter{section}{-1}

\section{\bf Introduction}\label{introduction}
\markboth{S. Azam, H. Yamaneh, M. Yousofzadeh}{Reflectable bases}

Over recent years there has been an increasing amount of
investigations on topics related to extended affine root systems,
extended affine Lie algebras and their generalizations. However,
in comparison with the inspiring models of finite and affine cases
(see \cite{Mac}, \cite{Hum}, \cite{K}, \cite{MP}), only a little
is known about the geometry of involved root systems (see
\cite{Hof2}). In the finite and affine Lie theory, the notion of
a ``root base" plays a crucial role in the study of not only the
corresponding geometry but also the whole theory, whereas it is
known that such a notion does not exist in general for the new
generalizations. In this work, we introduce the notions of a {\it
reflectable set} and a {\it reflectable base} for a {\it tame
irreducible  affine reflection system}, and we characterize
reflectable sets and reflectable bases for tame irreducible affine
reflection systems of reduced types, excluding types $E_{6,7,8}$.

Affine reflection systems are defined axiomatically in a way
similar to extended affine root systems. In fact, the axioms are
the same as those for  extended affine root systems, except that the underlying  finite dimensional vector space is replaced with an
arbitrary abelian group; see Definition \ref{newdef}. The notion
of an affine reflection system introduced here is more general
than the one defined by O. Loos and E. Neher \cite{LN2}, in fact
these two definitions are equivalent if the ground abelian group in our
definition is torsion free; see Remark \ref{clarify}. Affine
reflection systems include  extended affine
root systems \cite{AABGP}, locally extended affine root systems
\cite{MY}, and root systems extended by an abelian group
\cite{Y2}.

A reflectable set for a tame  irreducible  affine refection system
is a subset $\Pi$ of non-isotropic roots such that any
non-isotropic root can be obtained by reflecting a root of $\Pi$
relative to hyperplanes determined by $\Pi$. A reflectable base is
a minimal reflectable set; see Definition \ref{reflectable-base}.
One knows that any ``root base'' for a finite or affine root
system is a reflectable base, and it follows from our results in
this work that, for the types under consideration, the notion of a
reflectable base for a locally finite root system coincides with
the one introduced by Y. Yoshii \cite{Y3}. Also it is proved in
\cite{MS}, \cite{A2} and [LN1, Lemma 5.1] that any extended affine
root system or locally finite root system of reduced type
possesses a reflectable base. For extended affine root systems of
nullity $2$, a notion of ``root base'' is introduced by K. Saito
\cite{Sa}, and some of its algebraic and geometric features are
studied, however this notion of root base is not in general a
reflectable base as it might fail to have the minimality
condition. In the literature, one finds several other related terms
such as ``generalized bases'' \cite{NS}, ``integral bases''
\cite{LN1}, ``grid bases'' \cite{N2}, etc., each resembles, in
some aspects,  the usual notion of a root base for finite and
affine cases. It is shown in \cite{N2} that any grid base for a
3-graded root system is a reflectable base.

Reflectable bases have been appeared, though not necessarily with
such a name, in different contexts such as the description of root
systems, presentations of Weyl groups and  presentations of Lie
algebras \cite{Sa}, \cite{SaT}, \cite{AS1}, \cite{AS2},
\cite{AS3}, \cite{AS4}, \cite{Hof1} and \cite{SaY}. In the level of Lie algebras
and Weyl groups, an essential part of generators for a given
presentation is obtained from reflectable bases.  The Weyl group
acts naturally on the class of reflectable bases of an affine
reflection system. In the finite and affine cases, some orbits of
this action play a very important role, namely in the finite case
the class of root bases forms exactly one orbit of this action and
for the affine case, this class is the union of exactly two orbits; see \cite{K}. From this point of view, the study of orbits of
reflectable bases seems to be an interesting subject of research,
this has been one of our motivations for the study of such
objects. It is also known that the Weyl group of an affine
reflection system is not in general a Coxeter  group, for example
no extended affine Weyl group of nullity greater  than or equal 2
is a Coxeter group; see \cite{Hof1}. Nevertheless, non-Coxeter Weyl  groups are revealed to  have interesting geometric structures as well \cite{Hof2}. This has
been another motivation for us to do this work.

In this work, we give a full characterization of reflectable bases
for tame irreducible affine reflection systems of the  types under
consideration, namely types $A$, $B$, $C$, $D$, $F$ and $G;$ such
a characterization has  not been  known even for the finite and
affine cases. To establish  such a  characterization for   types
$E_{6,7,8},$ one needs further  investigations which  require an
independent study.

The paper is arranged as follows. In Section 1, we introduce
axiomatically the notion of an affine reflection system $R$ in the
same framework of an extended affine root system (see Definition
\ref{newdef}) and describe the structure of $R$ in the same lines
of \cite{AABGP} in terms of an involved locally finite root system
and some (pointed) reflection subspaces of the ground abelian
group; see Theorem \ref{structure-thm}. The relation between
affine reflection systems and certain other generalizations of
finite and affine root systems is clarified in Remark
\ref{clarify}. The notions of reflectable sets, reflectable bases
and integral bases are defined for affine reflection systems, and
some preliminary results are obtained for them. Some special
subsets of roots, which play an important role in our
characterization of reflectable sets and reflectable bases, are
introduced, including {\it coset spanning sets}, {\it coset bases}
and  {\it strong coset spanning sets}.

In Section 2, we characterize reflectable bases, and reflectable
sets for irreducible locally
 finite root systems of types under
consideration. For simply laced cases, a subset of nonzero roots
is a reflectable base if and only if it is a minimal set of
generators for the ground abelian group (see Proposition
\ref{supp1}) if and only if its image, under the canonical map, is
a basis for the $\bbbz_2$-vector space $\la\rd\ra/2\la\rd\ra,$ where $\la\cdot\ra$ denotes the $\bbbz-$span. This is what we call a coset basis; see Propositions \ref{a-new}
and \ref{asim1-new}. For non-simply laced types, the reflectable
sets and bases are characterized in terms of coset spanning sets
and coset bases of short and long roots in some appropriate vector
spaces. For example, according to our characterization, a subset
of an irreducible locally finite root system of type $B$, of rank
$\geq 3$, is a reflectable base if and only if it is the union of
a short root and a coset basis of long roots in $2\la \rd_{sh}\ra;$ see Proposition \ref{pro1-b>=3}. The section
contains also an expected but totally non-trivial fact that any
reflectable set contains a reflectable base; see Corollary
\ref{end5} and Propositions \ref{a-new} and \ref{asim1-new}. As a
by-product of our results in this section we see that any
reflectable base for an irreducible locally finite root system, of
types under consideration, is an integral base. In other words,
any reflectable base is a basis for the underlying free abelian
group. It is also shown that, for non-simply laced
 types, the cardinality $|\dot\Pi|$ of a reflectable base $\dot\Pi$ can be characterized in terms of dimensions of some vector spaces over Galois fields, namely
 $|\dot\Pi|=\dim(\la\rds\ra/\la\rdl\ra)+\dim(\la\rdl\ra/\la \rho\rds\ra),$
where here   $\rho$ stands for
the ratio of the long square root length to the short square root
length in $\rd$; see Remark \ref{cardinality-thm}.

In Section 3, using our results for locally finite root systems,
we give a full characterization of reflectable sets and
reflectable bases for tame irreducible affine reflection systems,
of the types under consideration, in terms of strong coset
spanning sets and minimal strong coset spanning sets. To see a
flavor of our results, let $R$ be an tame irreducible affine
reflection system with the sets of short roots  $\rsh$ and  long roots
$\rlg$. Then for type $A_1$, a subset $\Pi$ of $\r$ is a
reflectable set (resp. a reflectable base) for $R$ if and only if
$\Pi$ is a  strong coset spanning set (resp. a minimal strong coset
spanning set) for $R^\times$ in $2\la \r\ra$; see Theorem
\ref{finala1}. For type $B_2$, a subset $\Pi$ of $R$ satisfying
$\la\Pi\ra=\la R\ra$ is a reflectable set (resp. reflectable base)
for $\r$ if and only if $\Pi_{sh}:=\Pi\cap R_{sh}$ is a strong
coset spanning set (resp. a minimal strong coset spanning set) for
$\r_{sh}$ in $\la\r_{lg}\ra$ and $\Pi_{lg}:=\Pi\cap R_{lg}$ is a
strong coset spanning set (resp. a minimal strong coset spanning
set) for $\r_{lg}$ in $2\la\r_{sh}\ra$; see Theorem \ref{propro}.

We hope our characterization of reflectable sets and reflectable
bases offers a new perspective to the study of geometry of affine
reflection systems. The authors dedicate this work to the memory
of V. Shahsanaei who passed away in a fatal car accident at the
very early stage of this project, he was supposed to be one of the
authors.

Some parts of this work were completed during   ``QFT, String
Theory \& Mathematical Physics"  program at Kavli Institute for
Theoretical Physics (KITPC), China, held in July-August, 2010. The
authors would like to  thank the  KITPC for the warm hospitality.

\section{\bf Affine reflection systems}\label{preliminaries}
Throughout this work,   all vector spaces are considered over $\bbbq$
unless otherwise mentioned. For a vector space $\v,$ by $\v^\star,$
we mean the dual space of $\v.$ Also for  a subset $X$ of a group
$A,$  we denote by $\la X\ra,$ the subgroup of $A$ generated by $X$
and by $\Aut(A),$ we mean  the group of automorphisms of $A.$ In this
work  for a set $K$, $|K|$ denotes the cardinality of $K$ and  we
use the notation $\biguplus$, for  disjoint union. We recall  from
\cite{L} that a {\it pointed reflection subspace} of an additive
abelian group $A$ is a subset $X$ of $A$ satisfying one of the
following equivalent conditions:

- $0\in X$ and $X-2X\sub X$,

- $0\in X$ and $2X-X \sub X$,

- $2\la X\ra\sub X$ and $2\la X\ra-X\sub X$,

- $X$ is a union of cosets of $2\la X\ra$ in $X$, including the trivial coset $2\la X\ra$.

Also a {\it symmetric reflection subspace} of $A$ is a subset $X$ of $A$ satisfying one of the following equivalent conditions:

- $X-2X\sub X$,

- $X=-X$ and $2X+X\sub X$,

- $X=-X$ and $2X-X\sub X.$

Let $A$ be an abelian group, by a {\it symmetric form} on $A$, we
mean a symmetric bi-homomorphism $\fm:A\times A\longrightarrow
\bbbq$. This means that $\fm$ is a group homomorphism on each
component and is symmetric, namely $(\a,\a')=(\a',\a)$ and
$$(\a+\b,\a'+\b')=(\a,\a')+(\a,\b')+(\b,\a')+(\b,\b')$$
for all $\a,\a',\b,\b'\in A$. We set $A^0:=\{\a\in A\mid
(\a,A)=\{0\}\}$, and we call it the {\it radical} of $\fm$. We
also set $A^\times:=A\setminus A^0$, $\bar{A}:=A/A^0$ {and take
$\bar{\;}:A\rightarrow\bar{A}$  to be  the canonical epimorphism. For a subset $Y$ of $A,$ we denote by $\bar Y,$ the image of $Y$ under the map $\bar{\;}.$
We  note that  $\Ab$ is a torsion free group. The form $\fm$ is
called positive semidefinite (resp. positive definite) if
$(\a,\a)\geq 0$ (resp. $(\a,\a)>0$) for all $\a\in
A\setminus\{0\}$,} in this case one can see that
$$A^0=\{\a\in A\mid (\a,\a)=0\}.$$
For a subset $B$ of $A$, we set $B^\times:=B\setminus A^0$ and
$B^0:=B\cap A^0$. For $\a,\b\in A,$ if {$(\a,\a)\neq0$, we set
$(\b,\a^\vee):=2(\b,\a)/(\a,\a)$ and if  $(\a,\a)=0$,} we set
$(\b,\a^\vee):=0$. A subset $X$ of $A$ is called {\it connected}
if it cannot be written as a disjoint union of two nonempty
orthogonal subsets.  The form $\fm$ induces a unique form on
$\bar{A}$ by
$$(\ab,\bb):=(\a,\b)\qquad\hbox{for } \a,\b\in A.$$
This form is positive definite on $\bar{A}$.

Next suppose that  $R$ is a subset  of $A$ satisfying  $\la
R\ra=A$ and  $(\a,\b^\vee)\in\bbbz$  for all $\a,\b\in \rcross.$
For $\a\in\rcross,$ we take $w_\a\in\Aut(A)$ to be defined by
$w_\a(\b)=\b-(\b,\a^\vee)\a,$ $\b\in A,$ and
 {call it the {\it reflection based on
$\a.$}} We define the {\it Weyl group} $\w$ of $R$ to be the
subgroup of $\Aut(A)$ generated by $w_\a$, $\a\in\rcross$. In a
similar way, one defines $w_{\ab}\in\Aut(\bar{A})$ and $\bar{\w}$,
the subgroup of $\Aut(\bar{A})$ generated by $w_{\ab}$,
$\ab\in\bar{R}\setminus\{0\}.$ One can see that
\begin{equation}\label{conj}
sw_\a s^{-1}=w_{s\a};\quad s\in\w\andd\a\in R^\times.
\end{equation}
For a subset $\PP$ of $\rcross,$ we set
\begin{equation}\label{end8}\w_\PP:=\la w_\a\mid\a\in\PP\ra\andd
\w_\PP\PP:=\{w(\a)\mid w\in\w_\PP,\a\in\PP\}.\end{equation}
\begin{DEF}\label{newdef}
Let $A$ be an abelian group equipped with a nontrivial symmetric
positive semidefinite form $\fm$.
Let $R$ be a subset   of $A.$ The triple $(A,\fm,R),$ or $R$ if
there is no confusion, is called an {\it affine reflection system}
if it satisfies the following 3 axioms:

(R1) $R=-R$,

(R2) $\la R\ra=A$,

(R3) for $\a\in R^{\times}$ and $\b\in R$, there exist $d, u\in
{\mathbb Z}_{\geq 0}$ such that
$$(\b+\bbbz\a)\cap R=\{\b-d\a,\ldots,\b+u\a\}\andd
d-u=(\b,\a^\vee).
$$
Each element of $R$ is called a {\it root.} Elements of $R^\times$
(resp. $R^0$) are called  {\it non-isotropic  roots} (resp. {\it
isotropic roots}). The affine reflection system $R$ is called {\it
irreducible} if it satisfies

(R4) $\rcross$ is connected.

 Moreover, $R$ is called {\it tame} if

(R5) $R^0\sub\rcross-\rcross$ (elements of $R^0$ are non-isolated).

\noindent Finally $R$ is called {\it reduced} if it satisfies

(R6)$\;\a\in R^{\times} \Rightarrow 2\a\not\in R$.

An affine reflection system $(A,\fm,R)$ is called a {\it locally
finite root system} if $A^0=\{0\}.$
\end{DEF}

\begin{rem}\label{rem1}
 If $R$ is a finite root system in a Euclidean space $E,$ then $R$
is  a locally finite root system in $\la  R\ra.$ Therefore we
consider a finite root system as a locally finite root system.
\end{rem}

Let $(A,\fm,R)$ be
an affine reflection system. Here we describe the structure of $R$ in the same lines of
\cite{AABGP}  for extended affine root systems.
Since the form is nontrivial, we have $A\not=A^0$. Now this together
with (R1)-(R3), implies that
\begin{equation}\label{zero}0\in R.\end{equation}


\begin{lem}\label{non-deg}
If $(A,\fm,R)$ is an  affine reflection system, then the induced
form on $\v:=\bbbq\ot_\bbbz\bar A$ is non-degenerate. In
particular, $1\ot \bar R$ is locally finite in $\v,$ that is, any finite dimensional vector subspace of $\v$ intersects $1\ot \bar R$ in a finite set.
\end{lem}
\proof Denote the induced form on $\v$ by $\fm$ again. We recall
that this form satisfies
$$(r\ot \bar a,s\ot \bar b)=rs(\bar a,\bar b);\;r,s\in\bbbq,\;a,b\in A.$$

Suppose that $n$ is a positive integer, $\{p_i/q_i\mid 1\leq i\leq
n\}\sub\bbbq,$ $\{a_i\mid 1\leq i\leq n\}\sub A$ and
$\sum_{i=1}^n(p_i/q_i)\ot \bar a_i$ is an element of the radical
of the form on $\v.$ Take $q:=\Pi_{j=1}^nq_j$ and
$p'_i:=p_i\Pi_{i\neq j=1}^nq_j,$ $1\leq i\leq n,$ then
$$(1/q)\sum_{i=1}^n1\ot p'_i\bar a_i=(1/q)\sum_{i=1}^n p'_i\ot \bar
a_i=\sum_{i=1}^n (p'_i/q)\ot \bar a_i=\sum_{i=1}^n (p_i/q_i)\ot
\bar a_i$$ is an element of the radical of the form which in turn
implies that $\sum_{i=1}^n1\ot p'_i\bar a_i$ is an element of the
radical. Therefore for all $a\in A,$ $(
\sum_{i=1}^np'_i\bar a_i,\bar a)=(1\ot \sum_{i=1}^np'_i \bar
a_i,1\ot\bar a)=0.$ This means that $\sum_{i=1}^np'_i \bar a_i$ is
an element of the radical of the form on $\bar A.$ Now as the form
on $\bar A$ is non-degenerate, we are done. For the last
assertion,    using the same argument as in \cite[Lem.
I.2.6]{AABGP}, one can see that
\begin{equation*}
\label{bound} -4\leq 2(\b,\a)/(\a,\a)\leq 4;\;\; \a,\b\in
\rcross.\end{equation*} This together with the first part of the
proof and the same argument as in \cite[Pro. 3.7]{MY} completes
the proof. \qed

\begin{pro}\label{last1}   Suppose that  $(A,\fm,R)$ is a
locally finite root system, then  $\tilde R:=1\ot
R\sub\v:=\bbbq\ot_{\bbbz}A$ satisfies the following

(a) $0\in \tilde R,$ $\hbox{span}_\bbbq\tilde R= \v$ and $\tilde
R$ is locally finite,

(b) for every $\tilde\a\in \tilde R\setminus\{0\},$ there exists
$\check{\tilde\a}\in\v^\star$ such that
$\check{\tilde\a}(\tilde\a)=2$ and $\tilde R$ is invariant under the
reflection $ s_{\tilde \a}: \v\longrightarrow \v$ mapping $v\in \v$
to $v-\check{\tilde\a}( v)\tilde\a,$

(c) $\check{\tilde\a}(\tilde\b)\in\bbbz,$ for
$\tilde\a,\tilde\b\in \tilde R\setminus\{0\}.$

Conversely, if $\v$ is a vector space and $\tilde R\sub \v$
satisfies (a)-(c) above, then $\v$ is  equipped with a  symmetric
positive definite bilinear form $\fm$ invariant under
$s_{\tilde\a},$ $\tilde\a\in \tilde R\setminus\{0\}.$ Moreover
setting $A:=\la\tilde R\ra,$ we have that $(A,\fm_{\mid_{A\times
A}},\tilde R)$ is a locally finite root system.

\end{pro}

\proof  Suppose that $(A,\fm,R)$ is a locally finite root system.
 Since  $\la R\ra=A,$ we get
that
\begin{equation*}\hbox{span}_\bbbq \tilde R=\v.\end{equation*}
Next using Lemma \ref{non-deg}, we get that
\begin{equation*}
\parbox{4in}{\begin{center}$\tilde R$ is  locally finite in
$\v.$\end{center}}
\end{equation*}

Now using (R3), we get that for $\a,\b\in \rcross,$ $2(1\ot\b,1\ot
\a)/(1\ot\a,1\ot\a)\in\bbbz$ and that  the linear map $id\ot
w_\a:\v\longrightarrow \v$ preserves $\tilde R.$ Now if we set
$\tilde\a :=1\ot\a$ and define $\check{\tilde\a}\in\v^\star$  by
$v\mapsto 2(v,\tilde\a)/(\tilde\a,\tilde\a),$ then one gets that (a)-(c)
are fulfilled. Conversely, suppose that $\v$ is a vector space and $\tilde R\sub\v$ satisfies
(a)-(c), then by \cite[Thm. 4.2]{LN1}, there is a positive definite symmetric
bilinear form on $\v$ with desired property. Now using \cite[\S
4]{LN1}, the proof of  Proposition  3.15 of \cite{LN1} and general
facts on finite root systems, we are done. \qed

\begin{pro}\label{last2} Let $A$ be an abelian group and $R$ be a  subset of $A.$ Then there is a positive definite symmetric form $\fm $ on $A$ such that the
triple  $(A,\fm,R)$ is a locally finite root system if and only if
$A$ is a free abelian group of rank $dim(\v)$ where
$\v:=\bbbq\ot_\bbbz A$ and  $R$ satisfies

(a) $0\in R,$ $\la R\ra= A$ and $R$ is locally finite in the sense
that  any subgroup of $A$ of finite rank intersects $R$ in a
finite subset,

(b) for  every $\a\in \rcross,$ there exists
$\check\a\in\hbox{Hom}_\bbbz( A,\bbbq)$ such that
$\check{\a}(\a)=2$ and $R$ is invariant under the reflection $
s_{\a}: A\longrightarrow A$ mapping $a\in A$ to $a-\check{\a}(
a)\a,$

(c) $\check{\a}(\b)\in\bbbz,$ for $\a,\b\in R^\times.$ \end{pro}

\proof  Suppose that $(A,\fm,R)$ is a locally finite root system,
then  Proposition \ref{last1} together with  \cite[Lem.  5.1]{LN2}
implies that $1\ot R$ possesses a subset $B$ satisfying the
followings:

1) $B$ is a   basis for $\v,$

2) each element of $1\ot  R$ can be written as a $\bbbz$-linear
combination of elements of $B.$

This shows that $\la 1\ot  R\ra=\la B\ra$ is a free abelian group. Thic in turn together with the fact that $A$ is torsion free and
$B\sub 1\ot R$ implies that $A=\la  R\ra$ is a free abelian group of
desired rank. Next suppose that  $ X$ is a subgroup of $A$ of finite
rank, then there is a finite subset ${\mathcal{C}}$ of $B$ such that
$ 1\ot X\sub\la {\mathcal{C}}\ra. $ But $(1\ot R)\cap (1\ot X)\sub
(1\ot R)\cap\la {\mathcal{C}}\ra$ and so $(1\ot R)\cap (1\ot X)\sub
(1\ot R)\cap \hbox{span}_\bbbq( {\mathcal{C}}).$ Now since
$\u:=\hbox{span}_\bbbq( {\mathcal{C}})$ is a finite dimensional
subspace of $\v$ and by  Proposition \ref{last1},  $1\ot R$
is locally finite, we get that $(1\ot R)\cap\u$ is finite. Therefore  $(1\ot
R)\cap (1\ot X)$ and so $R\cap X$ is
finite. In other words, (a) is satisfied. Next to see that  (b) and (c) are fulfilled,  for $\a\in \rcross$
define
$$\check\a\in Hom_\bbbz(A,\bbbq);\;\; a\mapsto 2(\a,a)/(\a,\a)$$ an use  (R3).

Conversely,
suppose that $R$ is a subset of a free abelian group $A$ for which
(a)-(c) are satisfied. We first show that $1\ot R$ is  locally
finite in $\v=\bbbq \ot _\bbbz A.$ Take $\{a_i\mid i\in I\}$ to be
a basis for the free abelian group $A$ and suppose that $\u$ is a
finite dimensional subspace of $\v,$ then there is  a finite
subset $J$ of $I$ such that $\u\sub\hbox{span}_\bbbq \{1\ot
a_j\mid j\in J\}.$ Now it follows from the facts that $\{1\ot
a_i\mid i\in I\}$ and $\{a_i\mid i\in I\}$ are  bases for the
vector space $\v$ and the free abelian group $A$ respectively,
that if $1\ot \a\in (1\ot R)\cap\u,$ then $\a\in B\cap R$ in which
$ B$ is the subgroup of $ A$ generated by $\{ a_j\mid j\in J\}.$
Then  as $R$ is locally finite in $ A,$ we get that $B\cap R$ is
finite and so $(1\ot R)\cap U$ is finite, i.e., $1\ot R$ is
locally finite in $\v.$ Now for $\tilde\a:=1\ot \a,$ $\a\in
R^\times,$ define
$$\check{\tilde\a}:\v\longrightarrow \bbbq;\; q\ot a\mapsto q\check\a(a);\; q\in\bbbq,a\in A.$$
One can easily check that conditions  (a)-(c) in Proposition
\ref{last1} are satisfied for $\tilde R:=1\ot R\sub \v.$ Therefore
$\v$ is equipped with a positive definite symmetric bilinear form $\fm$ such
that $(1\ot A,\fm_{\mid_{(1\ot A)\times(1\ot A)}},1\ot R)$ is a
locally finite root system. Now we are done as the map from $A$ to
$\v$ mapping $a\in A$ to $1\ot a$ is an embedding.\qed

\begin{cor}\label{cor-last}
If $(A,\fm,R)$ is an affine reflection system, then $(\bar
R,\fm,\Ab)$ is a locally finite root system. In particular, if $R$
is irreducible, the induced form on  $\v:=\bbbq\ot_\bbbz \bar A $
is positive definite.
\end{cor}
\proof For the  first statement, using Lemma \ref{non-deg} and a
minor modification of the proof of  Proposition \ref{last2}, we
are done. The second assertion follows from Proposition
\ref{last1}, \cite[Thm. 4.2]{LN1} and the facts that the form
$\fm$ on $\bar A$ is positive definite and that this form is
invariant under $w_{\bar \a},$ $\a\in \rcross.$\qed\medskip

Now suppose that $\rdot$ is a locally finite root system in an
abelian group  $\dot A.$ Using Proposition \ref{last2}, one gets that $\dot A$ is
a free abelian group. We define  the {\it rank} of $\rdot$ to be the
rank of $\dot A.$ A subset $\dot S$ of $\rdot$  is said to be a {\it
subsystem} of $\rdot$ if it contains zero and $\dot w_{\da}(\db)\in
\dot S$ for $\da,\db\in \dot S\setminus\{0\}.$ 
Two locally finite root systems $(\rdot,\dot A)$
 and $(\dot S,\dot B)$ are  said  to be isomorphic if there is a group isomorphism
$f:\dot A\longrightarrow \dot B$ such that $f(\rdot)=\dot S.$

Suppose that $I$ is a nonempty  index set  and $\dot A:=\op_{i\in
I}\bbbz\ep_i$ is the  free  abelian group on   the
set $I.$ Define the  form $$\begin{array}{c}\fm:\dot A\times\dot A\longrightarrow\bbbq,\\
(\ep_i,\ep_j):=\d_{i,j}, \hbox{ for } i,j\in I.
\end{array}$$
This is a positive definite symmetric  form on $\dot A$.  Next define
\begin{equation}\label{locally-finite}
\begin{array}{l}
\dot A_I:=\{\ep_i-\ep_j\mid i,j\in I\},\\
D_I:=\dot A_I\cup\{\pm(\ep_i+\ep_j)\mid i,j\in I,\;i\neq j\},\\
B_I:=D_I\cup\{\pm\ep_i\mid i\in I\},\\
C_I:=D_I\cup\{\pm2\ep_i\mid i\in I\},\\
BC_I:=B_I\cup C_I.
\end{array}
\end{equation}
One can see that these are irreducible locally finite root systems
in their $\bbbz-$span's. Moreover, every irreducible locally finite
root system of infinite rank is isomorphic to one of these root
systems; see Proposition \ref{last1} and \cite[\S4.14 and
\S8]{LN1}. Now for an irreducible locally finite root system
$\rdot,$ define
$$\begin{array}{l}
\dot R_{sh}:=\{\da\in \dot R^\times\mid (\da,\da)\leq(\dot\b,\dot\b);\;\;\hbox{for all $\dot\b\in \dot R$} \},\\
\dot R_{ex}:=\rdot\cap2\dot R_{sh},\\
\dot R_{lg}:=\dot R^\times\setminus(\dot R_{sh}\cup\dot R_{ex}).
\end{array}$$
The elements of $\rdot_{sh}$ (resp. $\rdot_{lg},\rdot_{ex}$) are
called {\it short roots} (resp. {\it long roots, extra-long
roots}) of $\rdot$.  Using Proposition \ref{last1} and \cite[Cor.
5.6]{LN1}, one gets that each two nonzero roots of $\rdot $ of the
same length are conjugate under the Weyl group of $\rdot.$ In the
following, for not overusing the notations, irreducible finite
root systems of types $A_{\ell-1},$ $B_\ell,$ $C_\ell,$ $D_\ell$
and $BC_\ell$ will be denoted by $\dot A_I,$ $B_I,$ $C_I,$ $D_I$
and $BC_I$ respectively in which $I$ is an index set of
cardinality $\ell.$ We also refer to locally finite root systems
of types $\dot A_I,D_I,E_6,E_7$ and $E_8$ as simply laced types.

\begin{DEF}\label{type}
{\rm Considering Corollary \ref{cor-last}, for an affine
reflection system $R,$ we call the type and the  rank of $\bar R$
to be the {\it type} and the {\it rank} of $R$, respectively.}
\end{DEF}

The following proposition is  a generalization of  Proposition 5.9
of \cite{Hof1} to affine reflection systems.

\begin{pro}\label{end7}
Suppose that $A$ is an additive abelian group equipped with a
positive semidefinite symmetric form $\fm$ and $B$ is a subset of
$A^\times$ satisfying
\begin{itemize}
\item $\w_BB\sub B,$ \item $\bar B$ is an irreducible locally
finite root  system in $\hbox{span}_\bbbz \bar B.$
\end{itemize}
Then $R:= B\cup((B-B)\cap A^0)$ is a tame irreducible  affine
reflection system in $\hbox{span}_\bbbz(B).$
\end{pro}
\proof It is easy to see that (R1), (R2), (R4) and (R5) hold, so
we just need to check that (R3) is satisfied.  Let $\a\in
R^\times$ and $\b\in R$. We will look at three different cases.

\underline{\bf Case I:} $\ab$ and $\bb$ are linearly independent
over $\bbbz$, that is the subgroup of $\Ab$ generated by $\ab,\bb$
is of rank $2:$
Set $R_{\a,\b}:= R\cap(\bbbz\a\oplus\bbbz\b)$. Since the form on
$\bar A$ is positive definite and $\bar \a,\bar\b$ are
$\bbbz-$linearly independent, it is easy to see that the
restriction of the form to $\bbbz\a\op\bbbz\b$ is positive
definite. Now since  $(\b',{\a'}^\vee)\in\bbbz$ for all
$\a',\b'\in R_{\a,\b}$, it follows using the same argument as in
\cite[Lem. I.2.6]{AABGP} that $R_{\a,\b}$ is finite. Therefore it
is a finite root system in $\bbbz\a\oplus\bbbz\b$ and so (R3)
holds in this case.

\underline{\textbf{Case II:}} $\b\in R^0:$ If  $n\in\bbbz$ and
$\b+n\a\in R$,  then since $\w_BB\sub B,$ $\b-n\a=w_\a(\b+n\a)\in
B.$ So the set of integer numbers $n$ for which
$\bb+n\ab=n\ab\in\Rb$ is $\{0\}$, $\{0,\pm1\}$ or
$\{0,\pm1,\pm2\}$. But if $\b\pm 2\a\in R,$ again using the fact
that $\w_BB\sub B,$ one gets that $\b\pm\a=w_{\b\pm 2\a}(-\a)\in
B.$  So $\{\b+n\a\mid n\in\bbbz\}\cap R$ is equal $\{\b\}$,
$\{\b-\a,\b,\b+\a\}$ or $\{\b-2\a,\b-\a,\b,\b+\a,\b+2\a\}$.
Therefore  in each case (R3) holds.

\underline{\bf Case III:} $\ab$ and $\bb$ are nonzero and linearly
dependent over $\bbbz$:  Since $\bar R=\bar B$ is a locally finite
root system, we only need to consider three cases $\bar \a=\bar
\b$, $\bar\a=2\bar\b$ or $\bar\b=2\bar\a$. We show that in each
case the $\a$-string through $\b$ is of the desired form.  If
$\bar\a=2\bar\b,$ then  $w_\a(\b)=\b-\a$ and so the $\a$-string
through $\b$ is nothing but $\{\b-\a,\b\}$ as for $r\in\bbbz^{\leq
-2}\cup\bbbz^{\geq 1},$ $\overline{\b+r\a}=k\bar\b$ for some
integer $k$ with $|k|\geq 3$ which is not an element of $\bar R$.
For two other cases, we note that, by Case II, we are done if we
show that the mentioned string intersects $R^0$. We carry out this
as follows:

(i) $\bar\a=\bar\b:$ In this case  we have $\a-\b\in (B-B)\cap
A^0= R^0.$

(ii)  $\bar\b=2\bar\a:$ In this case the $\a-$string through $\b$
is a subset of $\{\b-4\a,\b-3\a,\b-2\a,\b-\a,\b\}.$ The following three cases can happen:
\begin{itemize}
\item $\b-2\a\in R:$ Since $\b-2\a\in A^0,$ $\b-2\a\in R^0$ and so
we are done.

\item $\gamma:=\b-3\a\in R:$ In this case $\bar \gamma=-\bar\a$
and so as in  Case III(i), we have $\b-2\a=\gamma+\a\in R^0.$

\item $\eta:=\b-4\a\in R:$ In this case, $\bar\eta=-2\bar\a.$ This
implies that $\eta\in B$ and so
$\gamma=\b-3\a=w_\eta\a\in\w_BB\sub B.$ Now using the previous
case, we are done.
\end{itemize}
This completes the proof.\qed
\begin{thm}\label{structure-thm}
Suppose that $(\dot A,\fm,\dot R)$ is an irreducible   locally
finite root system and $G$ is an abelian group.  If
$\rdl\not=\emptyset$, set
\begin{equation}\label{k}
\rho:=(\db,\db)/(\da,\da),\qquad (\da\in\rds,\;\db\in\rdl).
\end{equation} Let $S$, $L$, $E$ be subsets of $G$ satisfying
the conditions $(\star)$ and $(i)-(iv)$ below:
$$(\star)\;\;\;
\hbox{\parbox{3.5in}{$S,L$ (if $\rdl\neq \emptyset$) are pointed
reflection subspaces of $G$ and $E$ (if  $\rde\neq \emptyset$) is
a symmetric reflection subspace of $G$,}}$$
$\begin{array}{ll}(i)& \la S\ra =G,\vspace{2mm}\\
(ii)& S + L\subseteq S,\; L +\rho S \subseteq L\;
(\rdl\not=\emptyset),\vspace{2mm}\\
&S + E \subseteq S,\; E + 4S \subseteq E\;(\rd=BC_1),\vspace{2mm}\\
&L + E \subseteq L,\; E + 2L \subseteq E \;(\rd=BC_I,\;\;|I|\geq 2),\\
(iii) &\hbox{if  $\rdot$ is not  of  types
$A_1,$ $B_I$ or $BC_I,$ then $S=G,$}\vspace{2mm}
\\
(iv)& \hbox{if $\rdot$ is of types $B_I,$ $F_4,$ $G_2$ or $BC_I$
with $|I|\geq 3,$ then $L$ is a subgroup of $G.$}
\end{array}$
\newline\newline
Extend $\fm$ to a form on $A:=\dot A\op G$ such that
$(A,G)=(G,A)=\{0\}$, and set
$$(\star\star)\;\;\;R:=(S+S)\cup \underbrace{(\rds+S)}_{R_{sh}}\cup(\underbrace{\rdl+L}_{R_{lg}})\cup(\underbrace{\rde+E}_{R_{ex}})\sub
A$$where if $\rdl$ or $\rde$ is empty, the corresponding parts
vanish. Then $(R,\fm, A)$ is a tame irreducible  affine reflection
system. Conversely, suppose that $(R,\fm, A)$ is a tame irreducible
affine reflection system, then there is an irreducible locally
finite root system $\dot R$ and subsets $S,L,E\sub G:=A^0$ as in
$(\star)$ satisfying (i)-(iv) such that $R$ has an expression as
$(\star\star)$.
\end{thm}

\proof Using Proposition \ref{end7}, we get the first implication.
Conversely  suppose that $(A,\fm,R)$ is a tame irreducible affine
reflection system. Since $\Ab$ is torsion free, one can identify
$\bar\a\in \Ab$ with $1\ot \bar\a\in \bbbq\ot_\bbbz\Ab.$ So from
now on  we consider $\bar A$ as a subset of the $\bbbq$-vector
space $\v=\bbbq\ot_\bbbz\bar A.$ Now by Corollary \ref{cor-last},
$\bar{R}$ is a locally finite root system in $\v.$ Considering the
proof of Proposition \ref{last2} and using  Proposition
\ref{last1} and \cite[Lem. 5.1]{LN2}, one gets that $\Ab$ is a
free abelian group and there is a  basis $\bar\Pi\sub \bar R$ for
$\Ab$ satisfying $\w_{\bar{\Pi}}\bar{\Pi}=\Rb^\times$ (sets of
this form will be called integral reflectable bases later on). We
fix a pre-image $\dot{\Pi}\sub R$ of $\bar{\Pi}$ under the
projection map $^-$ and we set $\dot{A}:=\la\dot{\Pi}\ra$. Since
$\bar{\Pi}$ is a $\bbbz$-basis of $\Ab$, it follows that
$A=\dot{A}\oplus A^0$. We now set
$$\dot{R}:=(R+A^0)\cap \dot{A}.$$
Then an argument analogous to \cite[\S 2.2]{AABGP} shows that
$\dot{R}$ is a locally finite root system in $\dot{A}$  isomorphic
to $\Rb$. Again using exactly the same arguments as in \cite[\S
2.2]{AABGP}, one sees that if for $\da\in\dot{R}^\times,$ we set
$S_{\da}:=(\da+A^0)\cap R,$ then $S_{\dot\a}=S_{\dot\b}$ if $\dot\a$
and $\dot \b$ are of the same length and that  setting
$$\begin{array}{l}S:=S_{\dot\a};\;(\dot\a\in\rds),\\
L:=S_{\dot\a};\;(\dot\a\in\rdl\hbox{ if $\rdl$ is
nonempty}),\\
E:=S_{\dot\a};\;(\dot\a\in\rde\hbox{ if $\rde$ is
nonempty}),\end{array}$$  we have $(S,L,E)$ is as in $(\star)$ and
$R$ has an expression as in $(\star\star).$\qed
\bigskip

\begin{pro}\label{end9}
Suppose that $(A,\fm,R)$ is an affine reflection system and
$\mathcal{C}$ is a connected subset of $\rcross.$ Then

(i) $S:=\la \mathcal{C}\ra\cap R $ is an irreducible affine
reflection system in $\la \mathcal{C}\ra=\la S\ra.$

(ii) $T:=S^\times\cup((S^\times-S^\times)\cap A^0)$ is a tame affine
reflection system in $\la\mathcal{C}\ra=\la S^\times\ra=\la T\ra.$
\end{pro}
\proof (i) From the way $S$ is defined, we have
$\la\mathcal{C}\ra=\la S\ra.$ Now it is immediate that (R1) and
(R2) are satisfied. Next  suppose $\a\in S^\times$ and $\b\in S,$
then for $k\in\bbbz,$  $\b+k\a\in S$ if and only if $\b+k\a\in R.$
Now as (R3) holds  for $R,$ one gets (R3) for $S.$ Now we show that $S$ is irreducible. Suppose that
$S^\times=S_1\uplus\cdots\uplus S_t$ where $S_1,\ldots,S_t$ are
connected subsets of $S^\times$ with $(S_i,S_j)=\{0\}$ for $1\leq
i\neq j\leq t.$ Since $\mathcal{C}\sub S^\times$ is connected,
there is $1\leq j\leq t$ such that $\mathcal{C}\sub S_j.$ Now as
for each $1\leq k\leq t$ with $k\neq j,$ $(S_k,S_j)=\{0\},$ we
have $(S_k,\la \mathcal{C}\ra)=\{0\}.$ This implies that
$S^\times=\la\mathcal{C}\ra\cap R^\times\sub S_j$ which in turn
gives that $t=1$ and so $S^\times$ is connected.

(ii) One knows that $\w_{S^\times}S^\times\sub S^\times.$ Also as
$S$ is an affine reflection system in $\la S\ra,$ one gets using
Corollary \ref{cor-last} that the image of $S$ under the canonical
projection map $\la S\ra\longrightarrow \la S\ra/\la S\ra^0$ is an
irreducible locally finite root system in its $\bbbz-$span. Now By
Proposition \ref{end7}, $T:=S^\times\cup((S^\times-S^\times)\cap
A^0)=S^\times\cup((S^\times-S^\times)\cap \la S\ra^0)$ is a tame
affine reflection system in $\la S^\times\ra=\la T^\times\ra=\la
T\ra.$\qed

\begin{rem}\label{clarify}
(i) Let $\rd$ be an irreducible locally finite root system in an
abelian group  $\dot A$ and let $R=\{S_{\da}\}_{\da\in\rd^\times}$
be a root system extended by an abelian group $G$ in the sense of
\cite{Y2} where one replaces ``finite'' with  ``locally finite''.
Fix $\da\in\rds$, $\db\in\rdl$, if $\rdl\not=\emptyset$, and
$\dg\in\rde$, if $\rde\not=\emptyset$, and set $S:=S_{\da}$,
$L:=S_{\db}$ and $E:=S_{\dg}$. Then $S$, $L$ and $E$ satisfy all
conditions appearing in Theorem \ref{structure-thm}, and so the
set $R$ defined by ($\star\star$) is a tame irreducible affine
reflection system in $A:=\dot A\oplus G$. Conversely, let
$(A,\fm,R)$ be a tame irreducible  affine reflection system and
let $\rd$, $S$, $L$ and $E$ be as in the reverse part of Theorem
\ref{structure-thm}. For $\da\in\rd^\times,$ set $S_{\da}:=S$, if
$\da\in\rds$, $S_{\da}:=L$, if $\da\in\rdl$, and $S_{\da}:=E$, if
$\da\in\rde$. Then using Theorem \ref{structure-thm}, it is
straightforward to see that the collection
$$
R:=\{S_{\da}\}_{\da\in\rd^\times}$$ is a root system extended by
the abelian group $G:=A^0$.

(ii) Let $(A,\fm,R)$ be a tame irreducible  affine reflection
system and transfer the form from $A$ to the $\bbbq$-vector space
$\v:=\bbbq\otimes_\bbbz A$. Let $\v^0$ be the radical of the form
on $\v$ and set $\vb:=\v/\v^0$. Since $A=\dot A\op A^0,$ it is not
difficult to see that $\bbbq\otimes A^0$ can be identified with
the radical of the form on $\v$. It follows that we may naturally
identify $\vb$ with $\bbbq\otimes\bar A$. Using this
identification together with  Corollary \ref{cor-last}, we conclude that the
form on $\vb$ is positive definite. Therefore the form on $\v$ is
positive semidefinite.

Suppose now that $A$ is torsion free. Then one might identify $A$
with the subset $1\otimes A$ of $\v$, under the assignment
$a\mapsto 1\otimes a$, $a\in A$. Under this identification $R$ is
identified with $1\otimes R$. It is now easy to see that the
triple $(1\otimes R, \fm,\v)$ satisfies axioms (LR1)-(LR4) of
\cite[Def. 1.2]{AY}, and so is an irreducible locally extended
affine root system in the sense of \cite[Def. 1.2]{AY}. By
\cite[Pro. 1.3]{AY}, the set of non-isotropic roots of $(1\otimes
R,\fm,\v)$ is a locally extended affine root system in the sense
of \cite{Y3}. Conversely, if $(R,\fm,\v)$ is an irreducible
locally extended affine root system in the sense of \cite{AY},
then it follows immediately from definition that $R$ is a tame
irreducible affine reflection system in the torsion free abelian
group $A:=\la R\ra$.

(iii) In \cite{N1}, the author defines the term ``an affine
reflection system'' in a setting different from us. It follows
from \cite[3.11(b)]{N1} and \cite[Pro. 1.3]{AY} that a subset  $R$
of a vector space  is  a tame irreducible affine reflection system
in the sense of \cite{N1} if and only if it  is an irreducible
locally extended affine  root system in the sense of \cite{AY}.
\end{rem}

Here we introduce  some terminologies and recall  an elementary
fact about free abelian groups which will be used frequently in
the sequel. Suppose that $G$ is an additive abelian group and  $p$
is a prime number. For a subgroup $H$ of G with $pG\sub H\sub G,$
consider $G/H$ as a vector space over $\bbbz_p.$ Let  $K$ be a subset of  $G$ with $H\sub \la K\ra.$ We say a subset  $S$ of $ K$ is a {\it coset
spanning set} for $K$ in $H$ if $\{s+H\mid s\in S\}$ spans the
vector subspace $\la K\ra/H.$ The subset $S$ is called a {\it
coset basis} for $K$ in $H$ if $\{s+H\mid s\in S\}$ is a basis for
the vector subspace $\la K\ra/H.$ We also  call $\RR\sub K$ a {\it
strong coset spanning set} for $K$ in $H$ (with respect to $G$),
if $K=\cup_{x\in\RR}(x+H)\cap K.$ One can easily see that any
strong coset spanning set is a coset spanning set. Here is an
example of a coset spanning set which is not strong: Set
$G:=\bbbz\sg_1\oplus\bbbz\sg_2$, $H:=2G$ and $K:=G$. Then
$\RR:=\{\sg_1,\sg_2\}$ is a coset spanning set for $K$ in $2G$
which is not strong. In fact $\RR$ is a coset basis. The set
$\RR'=\{\sg_1,\sg_2,\sg_1+\sg_2\}$ is a minimal strong coset
spanning set for $K$ in $2G$.

\begin{lem}\label{final7}
 Suppose that $p$ is  a prime  number, $G$  a free abelian group and $H$  a subgroup
of $G$ satisfying $pG\sub H\sub G.$  Suppose that
$^-:G\longrightarrow G/pG$ and $[\cdot]:G\longrightarrow G/H$ are
canonical projection maps. Let $\{\a_i\mid i\in I\}\sub H$ and
$\{\b_j\mid j\in J\}\sub G$ be such that $\{\bar \a_i\mid i\in I\}$
is a linearly independent subset of the $\bbbz_p-$vector subspace
$H/pG$ of $G/pG$ and $\{[\b_j]\mid j\in J\}$ is a linearly
independent subset of  the quotient space $G/H.$ Then
$\{\a_i,\b_j\mid i\in I,j\in J\}$ is a $\bbbz-$linearly independent
subset of the free abelian group $G.$
\end{lem}

Let $(A,\fm,R)$ be a tame irreducible affine reflection system with
Weyl group $\w.$ As we have seen in Theorem \ref{structure-thm}, we
have $A=\dot{A}\oplus A^0$, where $\dot{A}$ is a free abelian
subgroup of  $A$ and $A^0$ is the radical of the form. Let
\begin{equation}\label{p}\mathfrak{p}:A\longrightarrow A^0\end{equation} be the canonical projection map from $A$
onto $A^0$.

For a subset $\PP$ of $\r^\times$, considering (\ref{end8}), we
set
$$\PP_{sh}:=\PP\cap\r_{sh}\andd\PP_{lg}:=\PP\cap\r_{lg}.$$

\begin{DEF} \label{reflectable-base} Let $\PP\subseteq\rcross.$

(i) We call ${\PP},$  a {\it reflectable set} for $R$, if $\w_\PP\PP=R^\times$.

(ii) The  subset $\PP$ is called a  {\it reflectable base}  if it
is a reflectable set and  no proper subset of  $\PP$ is a
reflectable set.

(iii) If $A$ is a free abelian group, the subset $\PP$ is called
an {\it integral base} for $R$ if $\PP$ is a basis for $A.$

(iv) 
%
If $A$ is  a free abelian group, we call $\PP$ an {\it integral
reflectable base} if $\PP$ is a reflectable base which is an
integral base as well.

\end{DEF}
We remark that if $\Pi$ is a nonempty subset of $\rcross$ and
$\a\in\Pi,$ then $\pm \a\in \w_\Pi\Pi$ and that  $\Pi$ is a
reflectable set (resp. reflectable base) if and only if
$(\Pi\setminus\{\a\})\cup\{-\a\}$ is a reflectable set (resp.
reflectable base). We refer to the latter property as the {\it sign
freeness} condition of reflectable sets (resp. reflectable bases).
We also remark that by Proposition \ref{last1} and \cite[Lemma 5.1]{LN2}, any reduced
locally finite root system possesses an integral  reflectable
base. In \cite{MS} and \cite{A2},  a reflectable base is
constructed for any  extended affine root system of reduced type.

\begin{lem}\label{connected}
Any reflectable set for $R$ is a connected generating set for $\la
R\ra$.
\end{lem}

\proof Let $\PP$ be a reflectable set. Clearly we have $\la R\ra=\la
R^\times\ra=\la\w_\PP\PP\ra\subseteq\la \PP\ra\subseteq\la R\ra$ and
so $\PP$ generates $\la R\ra$.

Next suppose that $\PP=\PP_1\cup\PP_2$ where $\PP_1,\PP_2$ are two
nonempty subsets with $(\PP_1,\PP_2)=\{0\}$. Then
$\w_\PP\PP_i\sub\la\PP_i\ra$ for $i=1,2$. So
$R^\times=\w_\PP\PP=\w_\PP\PP_1\cup\w_\PP\PP_2$ and
$(\w_\PP\PP_1,\w_\PP\PP_2)\subseteq
(\la\PP_1\ra,\la\PP_2\ra)=\{0\}$. This implies that $R^\times$ is
disconnected, a contradiction.\qed

\begin{lem}\label{21-11} Let $\Pi\subseteq\r^\times.$

(i) If  $\PP\sub\Pi$, $\b\in\Pi\setminus\PP$ and $\Pi'$ is  the set
obtained from $\Pi$ by replacing $\b$ with any element of the orbit
$\w_{\PP}\cdot\b$. Then $\w_{\Pi}\Pi=\w_{\Pi'}\Pi'$. In particular
$\Pi$ is a reflectable set if and only if $\Pi'$ is a reflectable
set.

(ii) If $\a,\a-\mathfrak{p}(\a)\in\Pi$, then $\Pi$ is a
reflectable set if and only if
$\Pi':=(\Pi\setminus\{\a\})\cup\{\a-2\mathfrak{p}(\a)\}$ is a
reflectable set.
\end{lem}

\proof (i) Take  $\a:=w'\b$ for some $w'\in\w_{\PP}$ and
$\Pi':=(\Pi\setminus\{\b\})\cup\{\a\}$. We show that
$\w_{\Pi'}=\w_\Pi$. For this,  it is enough to show that
$w_\gamma\in\w_{\Pi'}$ for any $\gamma\in\Pi$. If
$\gamma\not=\b$, then $\gamma\in\Pi'$ and so
$w_\gamma\in\w_{\Pi'}$. If $\gamma=\b$, then
$w_\gamma=(w')^{-1}w_{\a}w'\in\w_{\PP}\w_{\Pi'}\w_{\PP}\subseteq\w_{\Pi'}$.
Similarly, we have $\w_{\Pi'}\sub\w_\Pi.$ This completes the
proof.

(ii) Set $\PP=\Pi\setminus\{\a\}$. Then
$-\a+2\mathfrak{p}(\a)=w_{\a-\mathfrak{p}(\a)}(\a)\in\w_{\PP}\cdot\a$.
So by Part (i), with $\b=\a$, we get $\Pi$ is a reflectable set if
and only if $(\Pi\setminus\{\a\})\cup\{-\a+2\mathfrak{p}(\a)\}$ is
a reflectable set. But then by sign freeness, the latter is a
reflectable set if and only if $\Pi'$ is a reflectable set, as
required.\qed

\begin{lem}\label{general} Let $R$ be of one of  non-simply laced types and  $\a,\a_1,\ldots,\a_n\in\r^\times.$
Suppose that
$\{\a_{i_1},\ldots,\a_{i_t}\}=\{\a_1,\ldots,\a_n\}\cap\rsh$ and
$\{\a_{j_1},\ldots,\a_{j_m}\}=\{\a_1,\ldots,\a_n\}\cap\rlg$. Then
$$
w_{\a_1}\cdots w_{\a_n}(\a)\in\left\{
\begin{array}{ll}w_{\a_{i_1}}\cdots w_{\a_{i_t}}(\a)+\la\rlg\ra &\hbox{if }\a\in\rsh\\
w_{\a_{j_1}}\cdots w_{\a_{j_m}}(\a)+\rho\la\rsh\ra &\hbox{if }\a
\in\rlg.
\end{array}\right.
$$
\end{lem}

\proof Let $\a,\b\in\rlg$ and $\gamma\in\rsh$. We know that
$(\a,\gamma^\vee)\in \rho \bbbz.$ Therefore
$w_\gamma(\a)\in\a+\rho\la\rsh\ra.$ Then
$$
w_\gamma w_\b(\a)=w_\b(\a)+\rho\la\rsh\ra.$$ Now  it follows
using an  inductive process that  $w_{\a_1}\cdots w_{\a_n}(\a)\in
w_{\a_{j_1}}\cdots w_{\a_{j_m}}(\a)+\rho\la\rsh\ra$. If
$\a\in\rsh$, then an analogous argument as above, using the fact that
$(\b,\a^\vee)\in\bbbz$ if $\b\in\rlg,$ gives the other
implication.\qed

Note that if $R$ is of one of  non-simply laced types, it follows
from Theorem \ref{structure-thm} and the known facts on   locally
finite root systems (see Proposition \ref{last1} and \cite{LN1})
that
$$\rho\la\rsh\ra\sub\la\rlg\ra\subseteq\la\rsh\ra.
$$
Therefore,
\begin{equation}\label{vector-space}
\la\rsh\ra/\la\rlg\ra\hbox{ and }\la\rlg\ra/\rho\la\rsh\ra\hbox{
are two vector spaces over }\bbbz_\rho.
\end{equation}
\begin{lem}\label{general2} Let $R$ be of one of  non-simply laced types and $\Pi$ be a reflectable set for $R$.
Then $\Pi$ contains a subset $\MM$ such that
$\MM_{sh}=\MM\cap\rsh$ (resp. $\MM_{lg}=\MM\cap\rlg$)  is a coset
basis for $ \rsh$ in $\la\rlg\ra$ (resp. for  $\rlg$ in $\rho\la
\rsh\ra$).  Moreover, we have

(i)  $\{\mathfrak{p}(\a)\mid\a\in\MM^{lg}\}$ is a coset spanning
set for  $ L$ in $\rho\la S\ra,$

(ii) $\{\mathfrak{p}(\a)\mid\a\in\MM^{sh}\}$ is a coset spanning
set for
 $ S$ in $\la L\ra.$
\end{lem}

\proof From Lemma \ref{general}, it follows that, as subsets of
the abelian group $\la\r_{lg}\ra/\rho\la R_{sh}\ra$, we have
\begin{eqnarray*}
\la R_{lg}\ra/\rho\la R_{sh}\ra&=&
\la\w_{\Pi}\Pi_{lg}\ra/\rho\la R_{sh}\ra\\
&=&\big\la\{w\a+\rho\la R_{sh}\ra\mid w\in\w_\Pi,\a\in\Pi_{lg}\}\big\ra\\
&=& \big\la\{\a+\rho\la R_{sh}\ra\mid \a\in\Pi_{lg}\}\big\ra.
\end{eqnarray*}
Therefore $\Pi_{lg}$ is  a  {coset spanning set} for $\rlg$ in
$\rho\la \r_{sh}\ra.$ So $\Pi_{lg}$ contains a
 {coset basis $P_1$ for} $\rlg$ in $\rho\la
\r_{sh}\ra.$ An analogous argument shows that $\Pi_{sh}$ contains
a coset basis $P_2$ for $\rsh$ in $\la \rlg\ra.$ Now the set
$\MM:=P_1\cup P_2$ satisfies the first assertion of the statement
with $\MM_{sh}=P_1$ and $\MM_{lg}=P_2.$

Next, let $\d\in L$. Then
\begin{eqnarray*}
\d+\rho\la\rsh\ra\in\sum_{\a\in\MM_{lg}}\bbbz\a+\rho\la\rsh\ra&=&
\sum_{\a\in\MM_{lg}}\bbbz(\a-\mathfrak{p}(\a)+\mathfrak{p}(\a))+\rho\la\rds\ra+\rho\la
S\ra.
\end{eqnarray*}
From this it follows that
$\d\in\sum_{\a\in\MM_{lg}}(\mathfrak{p}(\a)+\rho\la S\ra)$,
showing that Part (i) holds. The argument for Part (ii) is
similar.\qed

\section{\bf Characterization of reflectable bases for  locally finite root systems}\setcounter{equation}{0}\label{subsub1-1}

In this section, we characterize reflectable bases for locally
finite root systems. This will be essential in obtaining the
characterization theorem for the general case. Let $\rd$ be a
locally finite root system of the form (\ref{locally-finite}). For
$\da=\sum_{i\in I}\eta_i\ep_i\in\rd$, we set
$$\supp(\da):=\{i\in I\mid\eta_i\not=0\}.$$ Also for a subset
$\dot{\mathcal{N}}$ of $\rd^\times,$ we set
$$\supp(\dot{\mathcal{N}}):=\cup_{\da\in\dot{\mathcal{N}}}\supp(\da).$$

\begin{lem}\label{a-d}
Let $\rd$ be an irreducible locally finite root system of type
$X=\dot{A}$,  $D$ or $B_I$ ($|I|\geq3$) and   $\PP$
be a coset spanning set for $\rd^\times$ in $2\la\rd\ra$ if
$X=\dot A,D$ and a coset spanning set for $\rdl$ in $ 2\la\rd_{sh}\ra$
if $X=B.$ Then $\supp(\PP)=I$ and $\PP$ is connected.
\end{lem}
\proof The first assertion in the statement   is easily checked. So we just need to show the second assertion. We carry out this in the following two
steps:

Step 1. If $\PP$ is disconnected, then there are two nonempty
subsets $\PP_1,\PP_2$ of $\PP$ with
$\supp(\PP_1)\cap\supp(\PP_2)=\emptyset$ such that
$\PP=\PP_1\cup\PP_2:$ Indeed, we claim that  there exists  either $\a\in
\PP_1$ such that $\supp(\a)\not=\supp(\b)$ for all $\b\in \PP_2$,
or  $\b\in \PP_2$ such that $\supp(\b)\not=\supp(\a)$
for all $\a\in \PP_1$. Suppose this does not hold, then for any
$\a\in \PP_1$, there exists a unique $\bar\a\in \PP_2$ such that
$\supp(\a)=\supp(\bar\a)$ and $(\a,\bar\a)=0$. Similarly, for any
$\b\in \PP_2$, there exists a unique element $\bar\b\in \PP_1$
such that $\supp(\b)=\supp(\bar\b)$ and $(\b,\bar\b)=0$. Now if
$\a_1,\a_2\in \PP_1$ and $\supp(\a_1)\cap\supp(\a_2)=\{i\}$ for
some $i\in I$, then $\bar\a_1,\bar\a_2\in \PP_2$ and
$\{i\}=\supp(\bar\a_1)\cap\supp(\bar\a_2)$, so
$(\a_1,\bar\a_2)\not=0$ which contradicts the fact that $(\PP_1,\PP_2)=\{0\}$.
This contradiction shows that for each $i\in I$, there exists a
unique $\a\in \PP_1$ such that $i\in\supp(\a)$, and in this case
$\bar\a$ is the unique element in $\PP_2$ such that $i\in\bar\a$.
Since the $\bbbz_2-$vector spaces $\la\rd\ra/2\la\rd\ra$ (if $X=A,D$)
and $\la\rdl\ra/2\la\rds\ra$ are of dimension greater than or
equal 2 and $\PP$ is a coset spanning set,   one finds
$\b\in \PP_1\setminus\{\a\}$. Let $\supp\{\a\}=\{i_\a,j_\a\}$ and
$\supp(\b)=\{i_\b,j_\b\}$. Then as we have seen above, $\a$ is the
unique element in $\PP_1$ having $i_\a$ or $j_\a$ in its support.
Similarly $\b$ is the unique element in $\PP_1$ having $i_\b$ or
$j_\b$ in its support. Also since $B$ is a coset spanning set for $\rd$
in $2\dot A$, we have $\ep_{i_\a}+\ep_{i_\b}\in\la \PP_1\ra+\la
\PP_2\ra+2\dot A$. Now define $f:\sum_{i\in
I}\bbbz\ep_i\rightarrow \bbbz_2$ by
$f(\ep_{i_\a})=f(\ep_{j_\a})=1$ and $f(\ep_i)=0$ for all $i\in
I\setminus\{i_\a,j_\a\}$. Then $f(\ep_{i_\a}+\ep_{i_\b})=1$ but
$f(\la \PP_1\ra)=\{0\}$, $f(\la \PP_2\ra)=\{0\}$ and $f(2\dot
A)=\{0\}$, a contradiction. This proves that our claim is true, so
without loss of generality, we may assume that
\begin{equation}\label{cure1} \hbox{there  exists
$\b\in \PP_2$ such that $\supp(\b)\not=\supp(\a)$ for all $\a\in
\PP_1$}.\end{equation} Now set
$$\begin{array}{c}\PP':=\{\b\in \PP_2\mid\supp(\b)=\supp(\a)\hbox{ for some }\a\in
\PP_1\},\\\PP'_1:=\PP_1\cup \PP'\andd\PP'_2:=\PP_2\setminus
\PP'.\end{array}$$ Then using (\ref{cure1}), we have  $$\PP=\PP'_1\cup \PP'_2,\; \PP'_1\not=\emptyset,\;\PP'_2\not=\emptyset,\;(\PP'_1,\PP'_2)=\{0\}\andd\supp(\PP'_1)\cup\supp(\PP'_2)=\emptyset.$$ So replacing $\PP_i$
with $\PP'_i$, $i=1,2$, if necessary, we may assume that
$\supp(\PP_1)\cap\supp(\PP_2)=\emptyset$.

Step 2. $\PP$ is connected: Suppose to the contrary that  $\PP$ is not
connected, then by Step 1, there are nonempty subsets
$\PP_1,\PP_2$ of $\PP$ with
$\supp(\PP_1)\cap\supp(\PP_2)=\emptyset$ such that
$\PP=\PP_1\cup\PP_2.$ Fix $u\in\supp(\PP_1)$ and $v\in\supp(\PP_2).$ Take $G:=\la\rd\ra$ if $X=A,D$ and  $G:=\la\rds\ra$ if
$X=B.$ Now define $f:\sum_{i\in I}\bbbz\ep_i\rightarrow \bbbz_2$
such that $f(\ep_i)=1$ for all $i\in \supp(\PP_1)$ and
$f(\ep_i)=0$ for all $i\in \supp(\PP_2)$. Since $\PP$ is a coset
spanning set, we have $\ep_u-\ep_v\in\la \PP_1\ra+\la \PP_2\ra+2G$
and so $f(\ep_u-\ep_v)=1$, $f(\la \PP_1\ra)=\{0\}$, $f(\la
\PP_2\ra)=\{0\}$ and $f(2G)=\{0\}$, a contradiction. This completes the proof.\qed

\medskip

\noindent \underline{\textbf{Type $B_I$:}} Suppose $\rd$ is a
locally finite root system of type $B_I.$

\begin{pro}\label{typeb2} Suppose that $|I|=2.$ A subset $\dot\Pi$
of $\rd^\times$ is a reflectable set if and only if
$\dot\Pi\cap\rds\neq\emptyset$ and $\dot\Pi\cap\rdl\neq\emptyset.$
Moreover a reflectable set $\dot\Pi$ is a reflectable base if and
only if $|\dot\Pi|=2.$
\end{pro}
\proof It is easy to check.\qed
\begin{pro}
\label{pro1-b>=3} Let  $|I|\geq3.$

(i) If $\dot\MM\sub\rd_{lg}$ is a  coset spanning set for $
\rd_{lg}$ in $\la 2 \rd_{sh}\ra$ and $i,j\in I$ with $i\not=j,$ then there
exist $w\in\w_{\dot\MM}$ and $\da\in\dot\MM$ such that
$\{i,j\}=\supp(w\da)$.

(ii)  Suppose that  $\dot\Pi$ is a reflectable set for $\rdot.$
Then $\dot\Pi_{lg}$ is a  coset spanning set for $\rd_{lg}$ in
$2\la\rd_{sh}\ra.$

(iii) Let $\dot\MM\subseteq\rdl$ be a  coset  spanning set for
$\rd_{lg}$ in $\la 2\rd_{sh}\ra$ and $\dot\a\in\rds.$ Then
$\dot\Pi:=\{\da\}\cup\dot\MM$ is a reflectable set for $\rdot.$
Moreover if $\dot\MM$ is a coset basis, then $\dot\Pi$ is a
reflectable base for $\rdot.$

(iv) Suppose that  $\dot\Pi$ is a reflectable base for $\rdot.$
Then $\dot\Pi$ contains  a short root $\da$ and  a   coset basis
$\dot\MM$ of $\rdl$ in $\la 2\rds\ra$ such that
$\dot\Pi=\{\da\}\cup\dot\MM.$

\end{pro}
\proof  
%
%
(i) 
%
%
%
By Lemma \ref{a-d}, $\dot\MM$ is connected and $\supp(\dot\MM)=I.$
 If $\dot\MM$ contains an element of support
$\{i,j\},$ there is nothing to prove. Also, if there exist $t\in
I$ and $\db,\dot\gamma\in\dot\MM$ with $\supp(\db)=\{i,t\}$ and
$\supp(\dot\gamma)=\{j,t\},$ then
$\{i,j\}=\supp(w_{\db}(\dot\gamma))$ and we are done. Otherwise,
since $\dot\MM$ is connected and $\supp(\dot\MM)=I$, we may find
 $\da_0,\ldots,\da_n\in\dot\MM$ such that $\supp(\da_0)=\{i,s_0\}$, $\supp(\da_n)=\{s_{n-1},j\}$ and
 $\supp(\da_k)=\{s_{k-1},s_k\}$ for $k=1,\ldots,n-1$. Then
 $\supp(w_{\da_0}\cdots w_{\da_{n-1}}(\da_n))=\{i,j\}$, as required.

(ii) See Lemma \ref{general2}.

(iii) Assume that $\da=\eta_0\ep_{i_0}$ for some $i_0\in I,$
$\eta_0\in\{\pm1\}$. By Part (i), for any $j\in
I\setminus\{i_0\},$ there are $\eta'_0,\eta_j\in\{\pm1\}$ such
that $\eta'_0\ep_{i_0}+\eta_j\ep_j\in\w_{\dot\MM}\dot\MM.$ Set
$\PP':=\{\ep_{i_0},\ep_{i_0}-\ep_{j}\mid j\in
I\}\sub\w_{\dot\Pi}\dot\Pi$. Then from \cite[Lemma 5.1]{LN2}, it
follows that
$$
\rdotcross=\w_{\PP'}\PP'\sub\w_{\dot\Pi}\dot\Pi\sub\rdotcross
$$
which shows that $\w_{\dot\Pi}\dot\Pi=\rdotcross,$ i.e., $\dot\Pi$
is a reflectable set. Now suppose $\dot\MM$ is  a coset basis. We
show that $\dot\Pi$ is a reflectable base. Suppose that
$\PP\sub\dot\Pi$ is a reflectable set, then by Part (ii),
$\PP_{lg}\sub\dot\MM$ is a  coset spanning set for $\rd_{lg}$ in
$2\la\rd_{sh}\ra$ and so by minimality of $\dot\MM$ (see Part
(i)), we get that $\PP_{lg}=\dot\Pi_{lg}=\dot\MM$ which in turn
implies that $\PP=\dot\Pi.$

(iv) Since $\dot\Pi$ is a reflectable base,
$\dot\Pi_{lg}$ is a coset spanning set for  $\rd_{lg}$ in
$2\la\rd_{sh}\ra$  by Part (ii). Fix $\da\in \dot\Pi_{sh},$ then by Part (iii),
$\PP:=\{\da\}\cup\dot\Pi_{lg}$ is a reflectable set. Now the
minimality of $\dot\Pi$ implies that $\dot\Pi_{sh}=\{\da\}.$ Next,
we claim that $\dot\Pi_{lg}$ is a coset basis for $\rd_{lg}$ in
$2\la\rd_{sh}\ra,$ indeed if it is not the case, then there is a
coset spanning set $\dot \PP\subsetneq\dot\Pi_{lg}$ for $\rd_{lg}$
in $2\la\rd_{sh}\ra.$ Now using Part (iii), $\{\da\}\cup\dot\PP$
is a reflectable set for $\rd.$ This  contradicts the minimality of
$\dot\Pi$ and so we are done.\qed

\vspace{1cm} \noindent \underline{\textbf{Types $\dot A_I
(|I|\geq2), D_I (|I|\geq 4):$}} Suppose that $\rd$ is a locally
finite root system of the types under consideration.  If $\rd$ is of
type $A_1,$ then $\rd^\times=\{\pm\da\}.$ In this case, $\{\da\},$
$\{-\da\}$ and $\{\pm\da\}$ are the only reflectable sets and
$\{\da\},$ $\{-\da\}$ are the only reflectable bases.  So from now
on we assume $\rd$ is a locally finite root system of type $\dot
A_I (|I|\geq3) $ or $ D_I (|I|\geq 4).$

\begin{pro}\label{supp1}

(a) Let $\dot\MM\sub\rd^\times$ be a generating set for $\la\rd\ra,$
then we have the followings:

(i) $\dot\MM$ is connected.

(ii) For  {$i_0,j_0\in I$ with $i_0\neq j_0$}, there exist
$\da\in\dot\MM$ and $w\in\w_{\dot\MM}$ such that
$\supp(w\da)=\{i_0,j_0\}$.

(iii) If $\rd$ is of type $D_I$, then there exist
$\da,\db\in\w_{\dot\MM}\dot\MM$ such that $\supp(\da)=\supp(\db)$
and $(\da,\db)=0$.

(iv)   $\dot\MM$   is a reflectable set.

(b) $\dot\Pi\subseteq\rd^\times$ is a reflectable base for $\rd$ if
and only if $\dot\Pi$ is   a minimal generating set for $\la\rd\ra.$
\end{pro}

\proof (a)(i) Since any subset of $\rd^\times$ generating
$\la\rd\ra$ is a spanning  set for $\rd^\times$ in $2\la\rd\ra,$
we are done using Lemma \ref{a-d}.

(a)(ii)  Fix  {$i_0,j_0\in I$ with $i_0\neq j_0.$} If $\dot\MM$
contains an element whose support is $\{i_0,j_0\},$ there is nothing
to prove, otherwise  set $\dot\MM(i):=\{j\in
I\mid\{i,j\}=\supp(\da)\hbox{ for some }\da\in\dot\MM\},$ $i\in I.$
Since $\dot\MM$ generates the abelian group $\la \rd\ra$, it follows
that $\dot\MM(i)\not=\emptyset$ for all $i\in I$. Now this together
with the connectedness of $\dot\MM$ (Part (a)(i)) implies that there
are $\da,\db\in\dot\MM$  and
$\da_0:=\da,\ldots,\da_{n}:=\db$ in $\dot\MM$ such that
$\supp(\da)\cap\dot\MM(i_0)\not=\emptyset$,
$\supp(\db)\cap\dot\MM(j_0)\not=\emptyset,$
$\supp(\da_0)=\{i_0,i\}$, $\supp(\da_n)=\{j,j_0\}$ and
$$\supp(\da_1)=\{i,j_1\},\;\supp(\da_2)=\{j_1,j_2\},\ldots,\supp(\da_{n-1})=\{j_{n-1},j\}.$$
Then $\{i_0,j_0\}=\supp\big(w_{\da_n}w_{\da_{n-1}}\cdots
w_{\da_1}(\da_0)\big).$  {This completes the proof.}

(a)(iii) Suppose to the contrary that there  not exist
$\dot\a,\dot\b\in\w_{\dot\MM}\dot\MM$ with
$\hbox{supp}(\da)=\hbox{supp}(\db)$ and $(\da,\db)=0.$
 {Then fix $i_0\in I.$ It follows
from Part (a)(ii) that  for any $i\in I\setminus\{i_0\},$ there is
$r_i\in\{\pm1\}$ such that $\ep_{i_0}+r_i\ep_i\in
\w_{\dot\MM}\dot\MM.$ This implies that for $i,j\in
I\setminus\{i_0\}$ with $i\neq j,$
$r_i\ep_i-r_j\ep_j\in\w_{\dot\MM}\dot\MM$ and so one gets that
$$\w_{\dot\MM}\dot\MM=\{\pm(\ep_{i_0}+r_i\ep_i)\mid i\in I\setminus\{i_0\}\}
\cup\{\pm(r_i\ep_i-r_j\ep_j)\mid i,j\in I\setminus\{i_0\},i\neq
j\}.$$} Now define $f:\sum_{i\in
I}\bbbz\ep_i\longrightarrow\bbbz_3$ by
$$\ep_{i_0}\mapsto1,\;\;\ep_i\mapsto2r_i;\;\;i\in I\setminus\{i_0\}.$$
 {Then } $f(\la
\w_{\dot\MM}\dot\MM\ra)=\{0\}$ which is a contradiction as
$\ep_{i_0}-r_i\ep_i\in\rdl\sub\la\dot\MM\ra\sub\la\w_{\dot\MM}\dot\MM\ra$
and $f(\ep_{i_0}-r_i\ep_i)=2.$ This completes the proof.

(a)(iv) Let $i,j\in I$ with $i\neq j$. By Part (a)(ii), there exist
$\da\in\dot\MM$ and $w\in\w_{\dot\MM}$ such that
$\supp(w\da)=\{i,j\}$. Now if we are in
 {case } $\dot{A}_I$, then the only
roots with support $\{i,j\}$ are $\pm(\ep_i-\ep_j)$, and so we are
done in this case.

Next suppose we are in  { case} $D_I$. By Part (a)(iii),
$\w_{\dot\MM}\dot\MM$ contains at least two roots of the forms
$\bd_1=\ep_i-\ep_j$ and $\bd_2=\ep_i+\ep_j$. Now for any $s\in
I\setminus\{i,j\}$, by Part (a)(ii), $\w_{\dot\MM}\dot\MM$
contains either $\ep_s-\ep_i$ or $\ep_s+\ep_i.$ But by acting
$w_{\db_1}$ and $w_{\db_2}$ on the one which belongs to
$\w_{\dot\MM}\dot\MM$, we get $\pm(\ep_s\pm
\ep_j)\in\w_{\dot\MM}\dot\MM$. Since $s\in I\setminus\{i,j\}$ was chosen
arbitrary, it follows that
$\pm(\ep_s\pm\ep_t)\in\w_{\dot\MM}\dot\MM$ for all $s,t\in I$ with
$s\neq t$.

(b) It follows immediately from Lemma \ref{connected} and Part
(a)(iv).\qed

\medskip

\begin{pro}\label{a-new} Let $\rd$ be an irreducible locally finite root system of type $\dot{A}_I$, $|I|\geq 2$, and $\dot\PP$ be a subset
of $\rd^\times$. Then $\dot\PP$ is a reflectable set (reflectable
base) for $\rd$ if and only if $\dot\PP$ is a coset spanning set
(coset basis) for $\rd^\times$ in $2\la\rd\ra$. In particular,  any
reflectable set contains a reflectable base. Moreover, any
reflectable base is an integral base.
\end{pro}

\proof Let $\dot\PP$ be a coset spanning set for $\rd^\times$ in
$2\la\rd\ra$. Then $\dot\PP$ contains a coset basis $B$ of
$\rd^\times$ in $2\la\rd\ra$. Clearly $\supp(B)=I$ and so by Lemma
\ref{a-d}, $B$ is connected. Then $S:=(\w_B B)\cup\{0\}$ is an
irreducible locally finite root system in $\la B\ra=\la S\ra$ of the
same type as $\rd$.
Now let $i_0,j_0\in I$ with $i_0\not=j_0.$ Since $\supp(S)=\supp(B)=I$, there exist
$\da,\db\in S$ with $i_0\in\supp(\da)$ and $j_0\in\supp(\db)$. Also since
$S$ is connected, we may find elements
$\da_0:=\da,\da_1,\ldots,\da_{n-1},\da_n:=\db$ in $S$ such that
$(\da_i,\da_j)\not=0$ if and only if $|i-j|\leq 1$.
Then $\pm w_{\da_n}\cdots w_{\da_1}\da_0=\pm(\ep_{i_0}-\ep_{j_0})$.
This shows that $\pm(\ep_{i_0}-\ep_{j_0})\in S$. Now since $i_0$ and
$j_0$ were chosen arbitrary, we get $S=\rd$. Thus $B$ and so
$\dot\PP$ is a reflectable set. Moreover, $B$ is a reflectable base
since otherwise it contains a proper subset $B'$ which is a
reflectable set. But then $\la B'\ra=\la B\ra=\la\rd\ra,$ and so
$B',$ which  is a coset spanning set for $\rd^\times$ in $2\la\rd\ra,$ is properly contained in the coset basis $B$, a contradiction.

Conversely, suppose $\dot\PP$ is a reflectable set (reflectable
base). Then $\la\dot\PP\ra=\la \rd\ra$ and so $\dot\PP$ is a coset
spanning set for $\rd^\times$ in $2\la\rd\ra$. If $\dot\PP$ is a
reflectable base and $\dot\PP$ is not a coset basis, then $\dot\PP$
contains a proper subset $B$ which is a coset basis for $\rd^\times$
in $2\la\rd\ra$. But by what we have seen above, $B$ is a
reflectable base. This contradics the minimality of $\dot\PP$.  The last assertion in the statement follows from the first assertion and Lemma \ref{final7}.\qed

%
%
%

\begin{lem}\label{typed}
Let $\rd$ be an irreducible locally finite root system of type  $D$, and $S$ be
a subsystem of $\rd$ of type $D$ such that $\supp(S)=\supp(\rd)$.
Then $S=\rd$.
\end{lem}

\proof By assumption, we may assume there exists an index set $K$
such that
$$S^\times=\{\pm(\ep'_t\pm\ep'_k)\mid t\not=k\in K\}
\sub \rd^\times=\{\pm(\ep_i\pm\ep_j)\mid i\not=j\in I\}.$$ Fix two
distinct $u,v\in K$ and set
$$B':=\{\ep'_{u}-\ep'_t\mid t\in K\setminus\{u\}\}\cup\{\ep'_u+\ep'_v\}.$$
Then one knows that $B'$ is a basis for the free abelian group $\la
S\ra$. Now let $t,t'\in K\setminus\{u\}$ and $t\not=t'$. Since
$S\sub \rd$, we have $\da:=\ep'_u-\ep'_t=r\ep_m+s\ep_n$ and
$\db:=\ep'_u-\ep'_{t'}=r'\ep_l+s'\ep_p$ for some
$r,s,r',s'\in\{\pm1\}$ and $m,n,l,p\in I$. Since $(\da,\db)\not=0$
and $\da,\db$ are $\bbbz$-linearly independent, one sees that the
set $\{m,n,l,p\}$ is of cardinality $3$. So, we may assume that
$\da=r\ep_m+s\ep_n$ and $\db=r'\ep_m+s'\ep_p$ with $p\not=n$. Now
since any other element of $B'$ of the form $\ep'_u-\ep'_t$ must be
non-orthogonal to both $\da$ and $\db$, we conclude that, we may
assume $K\sub I$, $u=m$ and
$$
B'':=\{\ep'_u-\ep'_t\mid t\in
K\setminus\{u\}\}=\{r_t\ep_u+s_t\ep_t\mid t\in K\setminus\{u\}\},$$
where $r_t, s_t\in\{\pm 1\}$. Next, consider $\dg:=\ep'_u+\ep'_v\in
S$. We know that $\dg$ is orthogonal to none  of elements of $B''$
except $\ep'_u-\ep'_v=r_v\ep_u+s_v\ep_v$. Therefore
$\dg\in\pm(r_v\ep_u-s_v\ep_v)$. So without loss of generality, we
may assume that
$$B'=\{r_t\ep_{u}+s_t\ep_t\mid t\in K\setminus\{u\}\}\cup\{r_v\ep_u-s_v\ep_v\}.$$
Therefore $K=\supp(B')=\supp(S)=\supp(\rd)=I$. Moreover for any
$t\in I=K$, we have
$$\pm(\ep_v\pm\ep_t)=w_{r_t\ep_{u}+s_t\ep_t}(\pm(\ep_u\pm\ep_v))\in S.$$
Since this holds for any $t\in I$, it follows again from the Weyl
group action that $\pm(\ep_i\pm\ep_j)\in S$ for all $i,j\in I$. Thus
$S=\rd$, as required.\qed

\begin{lem}\label{d3-new} Let $\rd$ be an irreducible locally finite root system of type $D$ and
$B$  a connected subset of $\rd^\times$ such that $\supp(B)=I$.
Then for any $i\in I$, the set $B\cup\{2\ep_i\}$ generates
$\la\rd\ra$.
\end{lem}

\proof Fix $i\in I$ and let $B':=B\cup\{2\ep_i\}$. Since $B$ is
connected and $\supp(B)=I$, for $j\in I$ there exist $\da_0,\da_n\in
B$ such that $i\in\supp(\da_0)$, $j\in\supp(\da_n)$ and a chain
$\da_0,\ldots,\da_n$ in $B$, connecting $\da_0$ to $\da_n$. Now
$2\ep_j\in\la\da_0,\ldots,\da_n,2\ep_i\ra\sub\la B'\ra$. Thus
$2\ep_j\in\la B'\ra$ for all $j\in I$. Next, using the connectedness
of $B$, for any $j,k\in I=\supp(B)$ with $j\neq k,$ we find
$\da,\db\in B$ with $j\in\supp(\da)$, $k\in\supp(\db)$  and a chain
in $B$ as above connecting $\da=\da_0$ to $\db=\db_n$. It follows
that either $\ep_j+\ep_k$ or $\ep_j-\ep_k$ belongs to the
$\bbbz$-span of this chain. Thus $\pm(\ep_j\pm\ep_k)\in\la B'\ra$
for all  $j,k\in I$, and so $\la B'\ra=\la\rd\ra$. \qed

\begin{pro}\label{asim1-new} Let $\rd$ be an irreducible locally finite root system of type $D_I$, $|I|\geq4$, and $\dot\PP$ be a subset
of $\rd^\times$. Then $\dot\PP$ is a reflectable set (reflectable
base) for $\rd$ if and only if $\dot\PP$ is a coset spanning set
(coset basis) for $\rd^\times$ in $2\la\rd\ra$. Moreover, any
reflectable set contains a reflectable base and any reflectable
base is an integral base.
\end{pro}

\proof Let $\dot\PP$ be a coset spanning set for $\rd^\times$ is
$2\la\rd\ra$. Then $\dot\PP$ contains a coset basis $B$ for
$\rd^\times$ in $2\la\rd\ra$. Clearly $\supp(B)=I$ and so by Lemma
\ref{a-d}, $B$ is connected. Then $S:=(\w_B B)\cup\{0\}$ is an
irreducible locally finite root system in $\la B\ra=\la S\ra$.
If $S$ is of type $D,$ then by Lemma \ref{typed}, $S=\rd$ and so
$B$ and consequently $\dot\PP$ is a reflectable set. Moreover, if
$\dot\PP$ is a coset basis, then as in the proof of Proposition
\ref{a-new}, we conclude that $\dot\PP$ is a reflectable base.

Now  considering the above argument, we are done if we show that $S$
can only be of type $D$. Suppose not, then the only possibility for $S$  is to  be of type $\dot{A}_J$ for some index set $J$. Suppose
this holds. If $\la S\ra=\la\rd\ra$, then by Proposition
\ref{supp1}, $S$ is a reflectable set for $\rd$, so
$S=\w_SS=\rd^\times$ which contradicts
the fact that $\rd$ is not of type $D$. So we have $\la S\ra\subsetneq\la\rd\ra$. Set $K:=\la S\ra$, $\bar{K}:=K/2K$, and for a subset $T$ of $K, $  denote
by $\bar{T}$ the image of  $T$  in $\bar{K}$, under
the canonical map. We know that $B$ is a coset basis for
$\rd^\times$ is $2\la \rd\ra$ and $K\subsetneq\la\rd\ra$,
therefore $\bar{B}$ is a basis for  the vector space $\bar{K}$.
Since $K\subsetneq \la\rd\ra$, we have $\bar{K}\subsetneq
\la\rd\ra/2K$. Now $\bar{B}\sub\bar{\dot\PP}$ and $\bar{\dot\PP}$
spans $\la\rd\ra/2K$, so $\bar{B}$ can be extended to a basis
$\bar{C}$ of $\la\rd\ra/2K$ such that $B\subsetneq C\sub\dot\PP$.
Let $\da:=r\ep_u+s\ep_v\in C\setminus B$, for some $r,s\in\{\pm
1\}$. Then $B\cup\{\da\}$ is $\bbbz$-linearly independent, in
particular $\da\not \in K=\la B\ra$. Now we note that $S$ is a connected
subsystem of $\rd$ and $\supp(S)=I.$ So  using a chain in $B$
connecting two elements whose supports contain $u$ and $v$, we
conclude that either $\da=r\ep_u+s\ep_v\in S\sub K$ or
$\db:=r\ep_u-s\ep_v\in S\sub K.$ This gives thst  $\db\in K$ as  $\da\not\in K.$
Now $\{\da\}\cup B$ is linearly independent and $\db\in K$, thus
$\{2\ep_u\}\cup B$ is $\bbbz$-linearly independent. Therefore by
Lemma \ref{d3-new}, $B\cup\{2\ep_u\}$ is a basis for $\la\rd\ra$.
So we may define a homomorphism $f:\la\rd\ra\rightarrow\bbbz_2$
such that $f(2\ep_u)=1$ and $f(B)=0$. Since $B$ is a coset basis
for $\rd^\times$ in $2\la\rd\ra$ we have $2\ep_u\in\la
B\ra+2\la\rd\ra$. But $f(2\ep_u)=1$ and $f(\la B\ra
+2\la\rd\ra)=\{0\}$, a contradiction.

Conversely, assume that $\dot\PP$ is a reflectable set
(reflectable base) for $\rd$. Then $\la\dot\PP\ra=\la \rd\ra$ and
so $\dot\PP$ is a coset spanning set for $\rd^\times$ in
$2\la\rd\ra$. If $\dot\PP$ is a reflectable base and $\dot\PP$ is
not a coset basis, then $\dot\PP$ contains a proper subset $B$
which is a coset basis for $\rd^\times$ in $2\la\rd\ra$. But by
what we have seen above, $B$ is a reflectable base. This contradicts
the minimality of $\dot\PP$.

The last assertion in the statement follows from our above
arguments together with Lemma \ref{final7}.\qed

\vspace{1cm} \noindent\underline{\textbf{Type $G_2$:}} Let $\rd$ be
a finite root system of type $G_2$.

\begin{pro}\label{G2}
A subset $\dot\Pi$ of $\rd^\times$ is a reflectable base for $\rd$
if and only if  $\dot\Pi=\{\dot\a,\dot\b\}$ where $\dot\a$ is a
short root, $\dot\b$ is a long root and $(\dot\a,\dot\b)\neq 0.$
\end{pro}

\proof Assume first that $\dot\Pi$ is a reflectable base for
$\rd$. From Lemma \ref{general2}, we know that $\dot\Pi$ contains
a subset $\dot\MM$ such that ${\dot\MM}_{sh}\sub\rds$ is a coset
basis for $\rds$ in $\la\rdl\ra$ and ${\dot\MM}_{lg}\sub\rdl$ is a
coset basis for $\rdl$ in $\la 3\rds\ra$. But one knows that both
$\la\rds\ra/\la \rdl\ra$ and $\la\rdl\ra/\la 3\rds\ra$ are one
dimensional $\bbbz_3$-vector spaces. So for any roots
$\da\in\dot\Pi_{sh}$ and $\db\in\dot\Pi_{lg}$, we may assume
$\dot\MM_{sh}=\{\da\}$ and $\dot\MM_{lg}=\{\db\}$. Since $\dot\Pi$
is connected, by Lemma \ref{connected}, it contains roots of
different lengths which are not orthogonal. So we may also assume
that $(\da,\db)\not=0$. Then the Coxeter graph associated to
$\dot\MM$ is a $G_2$ Coxeter graph and so $\dot\MM$ is a
reflectable set. By minimality of $\dot\Pi$, we get
$\dot\MM=\dot\Pi$.

Conversely, assume that $\dot\a\in\rds$ and $\dot\b\in \rdl$ are
such that $(\dot\a,\dot\b)\neq 0.$ Then  the Coxeter graph
associated to $\dot\Pi:=\{\dot\a,\dot\b\}$ is a $G_2$ Coxeter graph
and so $\dot\Pi$ is a reflectable set. Also if follows from the fact
that $\dot\Pi$ is of cardinality 2 that $\dot\Pi$ is a minimal
reflectable set, in other words, a reflectable base.\qed

 \vspace{1cm} \noindent \underline{\textbf{Type
$F_4$:}}  Let $\{\ep_1,\ldots,\ep_4\}$ be the standard basis for
$\bbbq^4$ and set $\rd:=\{0\}\cup\{\pm\ep_i,\pm(\ep_i\pm\ep_j)\mid
1\leq i\neq j\leq 4\}\cup\{(1/2)(r_1\ep_1+\cdots+r_4\ep_4)\mid
r_1,\ldots,r_4\in\{\pm1\}\}$ which is a  locally finite root
system of type $F_4.$

\begin{lem}\label{need}
(i) Let $r_1,r_2,r_3,r_4,s_1,s_2,s_3,s_4\in\{\pm1\}$ and set
$\dot\gamma_1:=(1/2)(r_1\ep_1+r_2\ep_2+r_3\ep_3+r_4\ep_4),\dot\gamma_2:=(1/2)(s_1\ep_1+s_2\ep_2+s_3\ep_3+s_4\ep_4),$
then $(\dot\gamma_1,\dot\gamma_2)=0$ if and only if
$|\{i\in\{1,2,3,4\}\mid r_i=s_i\}|=2.$ Moreover,
$\dot\gamma_1+\la\rdl\ra=\dot\gamma_2+\la\rdl\ra$ if and only if
either $(\dot\gamma_1,\dot\gamma_2)=0$ or
$\dot\gamma_1=\pm\dot\gamma_2.$

(ii)  If $\dot\gamma_1,\dot\gamma_2\in\rdl,$ then $\dot\gamma_1+2\la\rds\ra=\dot\gamma_2+2\la\rds\ra$
if and only if either $(\dot\gamma_1,\dot\gamma_2)=0$ or $\dot\gamma_1=\pm\dot\gamma_2.$

\end{lem}
\proof (i) It is an easy verification.

(ii) The implication $\Leftarrow$ is immediate.  To see the reverse
implication, let $\dg_1,\dg_2\in \rdl$ with $(\dg_1,\dg_2)\not=0$
and $\dg_1\not=\pm\dg_2$. Then $|\supp(\dg_1)\cap\supp(\dg_2)|=1$.
By symmetry on indices, we may assume $\dg_1=r\ep_2+s\ep_3$ and
$\dg_2=r'\ep_3+s'\ep_4$, where $r,r',s,s'\in\{\pm 1\}$. Now if
$\dg_1-\dg_2\in 2\la\rds\ra$, then we have
$r\ep_2+(s-r')\ep_3-s'\ep_4\in\bbbz(\ep_1+\cdots+\ep_4)+2\sum_{i=2}^4\bbbz\ep_i$,
which is absurd.\qed

\begin{pro}\label{f4-locally}
(i) Let $\dot\MM_1\sub\rdl$ be a  coset spanning  set for  $\rdl$ in
$2\la\rds\ra$ and $\dot\MM_2\sub\rds$ be a coset spanning set for
$\rds$ in $\la\rdl\ra,$ then $\dot\MM_1$ and $\dot\MM_2$ are
connected. Also $\dot\MM_2$ consists of at least two nonorthogonal
elements $\da_1,\da_2$ with $\hbox{supp}(\{\da_1,\da_2\})=4.$

(ii) If $\dot\MM_1\sub\rdl$ is a    coset spanning  set for $\rdl$
in $2\la\rds\ra$  and $\dot\MM_2\sub\rds$ is a  coset spanning set
for $\rds$ in $2\la\rdl\ra$ such that
$\PP:=\dot\MM_1\cup\dot\MM_2$ is connected, then $\PP$ is a
reflectable set for $\rd.$ Moreover, if $\dot\MM_1$ and $\dot\MM_2$
are coset bases  {for $\rdl$ in $2\la\rds\ra$ and for $\rds$ in
$\la\rdl\ra,$ respectively}, then $\PP$ is a reflectable base.
Furthermore, all reflectable sets and reflectable bases give rise
in this manner.
\end{pro}

\proof (i) We first note that the  $\bbbz_2$-vector spaces
$\la\rdl\ra/2\la \rds\ra$ and $\la\rds\ra/\la\rdl\ra$ are both of
dimension 2. Since  $\dot\MM_1$ is a  coset spanning set for
$\rd_{lg}$ in $2\la\rd_{sh}\ra,$ $\dot\MM_1$ contains a coset
basis $\mathcal{B}:=\{\db_1,\db_2\}.$ By Lemma \ref{need}(ii),
$(\db_1,\db_2)\neq0$ and $\db_1\neq\pm\db_2.$ This implies that
$\hbox{supp}(\{\db_1,\db_2\})$ is of cardinality 3. Now suppose
that $\db\in\dot\MM_1,$ then as $|\hbox{supp}(\dot\MM_1)|\leq4,$
there is $i\in\{1,2\}$ such that
$|\hbox{supp}(\db)\cap\hbox{supp}(\db_i)|=1.$ This means that any
element of $\dot\MM_1$ is connected to either $\db_1$ or  $\db_2.$ This
completes the proof of connectedness of $\dot\MM_1.$ Also as
$\dot\MM_2$ is a coset spanning  set for $\rd_{sh}$ in
$\la\rd_{lg}\ra,$ $\dot\MM_2$ contains a coset basis
$\{\da_1,\da_2\}$ and so  $|\hbox{supp}(\{\da_1,\da_2\})|=4.$
Without loss of generality, we assume that $\hbox{supp}(\da_1)$ is
of cardinality 4. We also mention that it follows from Lemma
\ref{need} that $(\da_1,\da_2)\neq 0.$
 {Now} we consider two cases, either
$|\hbox{supp}(\da_2)|=1$ or $|\hbox{supp}(\da_2)|=4.$ If
$|\hbox{supp}(\da_2)|=1,$ then for  any element $\da$ of
$\dot\MM_2,$ depending on $|\hbox{supp}(\da)|,$ $\dot\a$  is
connected to either $\da_1$ or   $\da_2$ and then the result follows.
So suppose that $|\hbox{supp}(\da_2)|=4,$ then any short root
$\da$ with $|\hbox{supp}(\da)|=1$ is connected to both
$\da_1,\da_2$ and any short root $\da$ with $|\hbox{supp}(\da)|=4$
is connected to either $\da_1$ or $\da_2,$ using Lemma
\ref{need}(i). This completes the proof.

(ii) We first note that as $\dot\MM_1$ is a coset spanning set for
$\rdl$ in $2\la \rds\ra$ and   $\la \rdl\ra/2\la \rds\ra$ is a
vector space of dimension 2,  $\dot\MM_1$ contains at least two
elements. Since $\PP$ and $\dot\MM_1\sub\rdl$ are connected, one
finds $\da_1\in\dot\MM_2$ and $\db_1,\db_2\in\dot\MM_1$ such that
$(\da_1,\db_1)\neq0$ and $(\db_1,\db_2)\neq0.$  Without loss of
generality, we suppose that  there are $r,s,m,n\in\{\pm1\}$ such
that $\db_1=r\ep_1+s\ep_2$ and $\db_2=m\ep_1+n\ep_3.$ Next take
$\dot\a_2\in\dot\MM_2$ to be such that $\{\da_1,\da_2\}$ is a
coset basis for $\rds$ in $\la\rdl\ra.$ We consider two cases
$|\hbox{supp}(\da_1)|=1$  and $|\hbox{supp}(\da_1)|=4.$ In the
former case,  take $\dot S:=\{\pm\ep_i,\pm(\ep_i\pm\ep_j)\mid
1\leq i,j\leq 3,i\neq j\}$ which is a subsystem of $\rd$ of type
$B_3$ and note that  the Coxeter graph associated to the nods of
either $\{\da_1,\db_1,\db_2\}$ or
$\{\da_1,\db_2,w_{\db_2}(\db_1)\}$ is the same as the Coxeter
graph  of type $B_3,$ then it follows that
$\PP':=\{\da_1\}\cup\{\db_1,\db_2\}$ is a reflectable set for
$\dot S.$  Therefore  $\dot S=\w_{\PP'}{\PP'}\sub\w_{\PP}\PP.$
Next, we note that as $|\hbox{supp}(\dot\a_1)|=1,$ we get that
$|\hbox{supp}(\dot\a_2)|=4$ and so
$\da_2=(1/2)(r_1\ep_1+r_2\ep_2+r_3\ep_3+r_4\ep_4)$ for some
$r_1,\ldots,r_4\in\{\pm1\}.$ Now as $-r_1\ep_1-r_2\ep_2\in\dot
S\sub\w_\PP\PP, $ we get that
$\db_3:=r_3\ep_3+r_4\ep_4=w_{\da_2}(-r_1\ep_1-
r_2\ep_2)\in\w_{\PP}{\PP}.$ Now again  the Coxeter graph
associated to the nods of either {$\{\da_1,\db_1,\db_2,\db_3\}$ or
$\{\da_1,\db_2,w_{\db_2}(\db_1),w_{w_{\db_2}(\db_1)}(\db_3)\}$} is
the same as the Coxeter graph of type $B_4.$ This in turn together
with Lemma \ref{21-11} implies that $\PP'':=\PP'\cup\{\db_3\}$ is
a reflectable set for the subsystem
$\rd_B:=\{\pm\ep_i,\pm(\ep_i\pm\ep_j)\mid 1\leq i,j\leq 4,i\neq
j\}.$ This means that $\rd_B=\w_{\PP''}\PP''\sub\w_\PP\PP.$ Now to
complete the proof in this case, we must show any element of
$\rd\setminus\{\da_2\}$ whose  support is of cardinality $4,$
belongs to $\w_\PP\PP.$ So take $s_1,\ldots,s_4\in\{\pm1\}$ to be
such that
$\db:=(1/2)(s_1\ep_1+\cdots+s_4\ep_4)\in\rd\setminus\{\da_2\}.$
Suppose $1\leq t\leq 4$ and $1\leq i_1\leq\cdots\leq i_t\leq 4$
are such that $s_{i_j}=-r_{i_j}$ for $1\leq j\leq t,$ then
$\db=w_{\ep_{i_1}}\ldots w_{\ep_{i_t}}\da_2\in\w_{\PP}\PP.$ This
completes the proof in the former case. In the latter case, we
first note that by Lemma \ref{need}, $(\da_1,\da_2)\neq0.$  Also
by Lemma \ref{21-11}, $\PP$ is a reflectable set if and only if
$(\PP\setminus\{\da_2\})\cup\{w_{\da_1}\da_2\}$ is a reflectable
set, so without loss of generality we assume $\hbox{supp}(\da_2)$
is of cardinality 1. Since $(\db_1,\da_1)\neq 0,$ we get
$\da_1=\pm(1/2)(r\ep_1+s\ep_2+r'\ep_3+s'\ep_4)$ for some
$r',s'\in\{\pm1\}.$ Set $\db_3:=w_{\da_1}(\db_1),$ then
$\hbox{supp}(\db_3)=\{3,4\}$ and so by Proposition
\ref{pro1-b>=3}, $\PP_1:=\{\da_2,\db_1,\db_2,\db_3\}$ is a
reflectable set for $\rd_B$ i.e.,
$\{\pm\ep_i,\pm(\ep_i\pm\ep_j)\mid1\leq i\neq
j\leq4\}=\w_{\PP_1}\PP_1$  and so the same argument as above
completes the proof.

Now suppose $\dot\MM_1$ and $\dot\MM_2$ are coset bases for $\rdl$
in $2\la\rds\ra$ and for $\rds$ in $\la\rdl\ra,$ respectively, but
$\PP=\dot\MM_1\cup\dot\MM_2$ is not a reflectable base. So  there
is $\dot\PP\sub\PP$ such that $\dot\PP$ is a reflectable set for
$\rd$.  Thus by Lemma \ref{general2}, $\dot\PP_{sh}$ and
$\dot\PP_{lg}$ are coset spanning sets for $\rd_{lg}$ in
$2\la\rds\ra$ and for $\rds$ in $\la\rdl\ra,$ respectively. But one
knows that $\dot\PP_{lg}\sub\dot\MM_1$ and
$\dot\PP_{sh}\sub\dot\MM_2,$ therefore we get
$\dot\MM_1=\dot\PP_{lg}$ and $\dot\MM_2=\dot\PP_{sh}$ and so
$\PP=\dot\PP.$ It is easy to see, using Lemma \ref{general2}, that any
reflectable set is connected and that it  is of the form
$\dot\MM_1\cup\dot\MM_2$ where $\dot\MM_1$ is a coset spanning set
for $\rdl$ in $2\la\rds\ra$ and $\dot\MM_2$ is a coset spanning
set for $\rds$ in $\la\rdl\ra.$ Now suppose that $\dot\Pi$ is a
reflectable base for $\rd.$ We have already seen that
$\dot\Pi=\dot\Pi_{lg}\uplus\dot\Pi_{sh},$ and that $\dot\Pi_{lg}$
is a coset spanning set for $\rdl$ in $2\la\rds\ra$ and
$\dot\Pi_{sh}$ is a coset spanning set for $\rds$ in $\la\rdl\ra.$
Since $\dot\Pi$ is connected, there are $\da\in\dot\Pi_{lg}$ and
$\db\in\dot\Pi_{sh}$ such that $(\da,\db)\neq 0.$ Take
$\dg\in\dot\Pi_{lg}$ and $\dot\eta\in\dot\Pi_{sh}$ such that
$\{\da,\dot\gamma\}$ is a coset basis for $\rdl$ in $2\la\rds\ra$
and $\{\db,\dot\eta\}$ is a coset basis for $\rds$ in
$\la\rdl\ra.$ Therefore, $\{\da,\db,\dg,\dot\eta\}\sub\dot\Pi$ is
a reflectable base for $\rd$ and so
$\dot\Pi=\{\db,\dot\eta,\da,\dg\}.$ This completes the proof.\qed

\vspace{5mm} \noindent \underline{\textbf{Type $C_I(|I|\geq3)$:}}
Let $\rd$ be a locally finite root system of type $C_I$ ($|I|\geq
3$) and $\rd^\vee$ be its dual root system, namely
$\rd^\vee=\{0\}\cup\{\da^\vee:=2\da/(\da,\da)\mid\da\in\rd^\times\}$.
It is known that $\rd^\vee$ is a locally finite root system of type
$B_I$.

\begin{pro}\label{c}
Let $\dot\Pi\sub\rd^\times$. Then $\dot\Pi$ is a reflectable base
for $\rd$ if and only if $\dot\Pi=\{\db\}\cup\dot\MM,$ where $\db$
is a long root and $\dot\MM$ is a  coset basis  for $\rds$ in $\la
\rdl\ra$.
\end{pro}

\proof Set $\dot\Pi^\vee:=\{\da^\vee\mid\da\in\dot\Pi\}$. Let
$\da\in\w_{\dot\Pi}\dot\Pi$, say $\da=s_{\da_1}\cdots
s_{\da_n}(\da_{n+1})$ where $\da_i$'s are in $\dot\Pi$. Then
$\da^\vee=s_{\da_1^\vee}\cdots s_{\da_n^\vee}(\da_{n+1}^\vee)$.
This shows that
$(\w_{\dot\Pi}\dot\Pi)^\vee=\w_{\dot\Pi^\vee}(\dot\Pi^\vee)$. From
this it follows that $\dot\Pi$ is a reflectable base for $\rd$ if
and only if $\dot\Pi^\vee$ is a reflectable base for $\rd^\vee$.
Moreover, it is easy to see that dual of a coset basis  for $\rdl$
in $\la 2\rds\ra$ is a minimal coset spanning set for $\rds$ in
$\la\rdl\ra$ and vice versa. Now the result follows immediately
from Proposition \ref{pro1-b>=3}.\qed

\begin{cor}\label{end5}
Suppose that $\rd$ is a locally finite root system of a non-simply
laced type. Then any reflectable set for $\rd$ contains a
reflectable base for $\rd$ and any reflectable base is an integral
base.
\end{cor}

\proof Using Lemma \ref{general2}, Propositions \ref{typeb2},
\ref{pro1-b>=3}, \ref{G2}, \ref{f4-locally} and \ref{c}, we get the
first assertion. For the second assertion, using Propositions
\ref{pro1-b>=3}, \ref{f4-locally}, \ref{G2} and \ref{c}, we get that
$\dot\Pi_{sh},\dot\Pi_{lg}$ are coset bases for $\rd_{sh}$ in
$\la\rd_{lg}\ra$ and for $\rdl$ in $\la\rho\rds\ra,$ respectively.
Now as $\rho\la\rds\ra\sub\la\rdl\ra\sub\la\rds\ra,$ we are done
using Lemma \ref{final7}. \qed

\begin{rem}\label{cardinality-thm}
Suppose that  $\rd$ is  a locally finite root system of a
non-simply laced type. By Corollary \ref{end5}, we know that any
reflectable base is an integral base. This in particular implies
that $|\dot\Pi|=\rank(\rd)$. Moreover, it follows from  the proof
of Propositions \ref{typeb2}, \ref{pro1-b>=3}, \ref{G2} and
\ref{f4-locally}   that $|\dot\Pi|$ can be characterized as
\begin{eqnarray*}|\dot\Pi|&=&\dim(\la\rds\ra/\la\rdl\ra)+\dim(\la\rdl\ra/\la
\rho\rds\ra)\\&=& \left\{\begin{array}{ll}
\hspace{-2mm}1+\dim(\la\rdl\ra/\la 2\rds\ra)&\hbox{if }X=B_I,\\
\hspace{-2mm}\dim(\la\rds\ra/\la\rdl\ra)+1&\hbox{if }X=C_I,\\
\hspace{-2mm}4&\hbox{if }X=F_4,\\
\hspace{-2mm}2&\hbox{if }X=G_2.
\end{array}
\right.\end{eqnarray*}
\end{rem}

%
\bigskip

We conclude this section with the following  two propositions
which  will not be used in the sequel but reveal some
interesting relations  between root systems of types $\dot A_I,$
$B_I$ and $D_I.$ Let $\rd$  be a locally finite root system of
type $B_I,$ we recall that $\rd$ is of the form
$\{0,\pm\ep_i,\pm(\ep_i\pm\ep_j)\mid i\neq j\in I\}.$ Take
$S_A:=\{\pm(\ep_i-\ep_j)\mid i,j\in I,i\neq j\}\cup\{0\}$ and
$S_D:=\{\pm(\ep_i\pm\ep_j)\mid i,j\in I,\;i\neq j)\}\cup\{0\}.$
Then $S_A$ is a locally finite root system of type $A_{I}$ and
$S_D$ is a locally finite root system of type $D_I$ (if $|I|\geq
4$).
\medskip

\begin{pro}\label{att-to-a-2}
(i) Let $\rd$ be a locally finite root system of type $B_I$ and
$S_A$ be the subsystem of  $\rd$ of  type $A_I$ as above. Suppose
that $\dot\PP\sub S_A^\times$ and fix  $i\in I.$   Then $\dot\PP$ is
a reflectable set for $S_A$ if and only  if
$\dot\Pi:=\dot\PP\cup\{\ep_i\}$ is a reflectable set for  $\rd.$
Moreover, $\dot\PP$ is a reflectable base for $S_A$ if and only if
$\dot\Pi$ is a reflectable base for $\rd.$
\end{pro}

\proof (i) Suppose that $\dot\PP$ is a reflectable set for $S_A$.
Then $\dot\PP$ is connected and $\la\dot\PP\ra=\la S_A\ra.$ So
$\dot\PP$ is a connected coset spanning set for
$\rdl$ in $\la 2\rds\ra$. Therefore by Proposition
\ref{pro1-b>=3}(iii), $\dot\Pi$ is a reflectable set for $\rd$.
Next suppose that $\dot\Pi$ is a reflectable set for $\rd.$ We
show that $S_A^\times=\w_{\dot\PP}\dot\PP.$ Suppose that $r,s\in
I$ with $r\neq s.$ We  prove  that $\ep_r-\ep_s\in
W_{\dot\PP}\dot\PP.$ Since $\dot\Pi$ is a reflectable set for
$\rd,$ there are $\da_1,\ldots,\da_n\in\dot\Pi$  and
$\dot\a\in\dot\PP$ such that $w_{\da_1}\ldots
w_{\da_n}(\da)=\ep_r-\ep_s,$  with  $n$  as small as possible (we refer to $w_{\da_1}\ldots w_{\da_n}(\da)$ as a ``reduce expression''). We
show that $w_{\da_1}\ldots w_{\da_n}(\da)\in\w_{\dot\PP}\dot\PP.$
Let $f:\dot A\rightarrow\bbbz$ be the homomorphism induced by the
assignment $\ep_j\mapsto 1$ for all $j\in I$. Now
$0=f(\ep_r-\ep_s)=f(w_{\da_1}\cdots w_{\da_n}(\da))$. Note that if
$\da_k=\ep_i$ for some $1\leq k\leq n$, then $(\da_k,
w_{\da_{k+1}}\cdots w_{\da_n}(\da))\not=0$ as the expression is
reduced and $f(\da_j)=0$ for $\da_j\in\dot\PP$. Thus
$0=f(w_{\da_1}\cdots w_{\da_n}(\da))=\sum_{t=1}^{p}k_t$, where
$k_t\in\{\pm2\}$ and $p$ is the number of $j$'s for which
$\da_j=\ep_i$. Therefore  $p$ is even and so without loss of
generality, we may assume that $\da_1=\da_n=\ep_i,$ and $\da_j\neq
\ep_i$ for $2\leq j\leq n-1.$ Now as the expression is reduced, we
have $(\ep_i,\ep_r-\ep_s)=(\da,\ep_r-\ep_s)\not=0$ and
$(\ep_i,\da)=(\dot\a_n,\da)\not=0$. Now one can see that
$w_{\da_1}w_{\da_2}\ldots w_{\da_n}(\da),w_{\da_2}\ldots
w_{\da_{n-1}}(\da)\in S_A^\times$ and
$\hbox{supp}(w_{\da_1}w_{\da_2}\ldots
w_{\da_n}(\da))=\hbox{supp}(w_{\da_2}\ldots w_{\da_{n-1}}(\da))$.
Now the minimality of $n$ implies that $w_{\da_1}w_{\da_2}\ldots
w_{\da_n}(\da)=-w_{\da_2}\ldots w_{\da_{n-1}}(\da).$ Setting
$\db_1:=w_{\da_2}\ldots w_{\da_{n-1}}(\da),$ we have
$w_{\da_1}w_{\da_2}\ldots w_{\da_n}(\da)=w_{\db}w_{\da_2}\ldots
w_{\da_{n-1}}(\da)\in\w_{\dot\PP}\dot\PP.$ This completes the
proof of the first assertion. The second assertion  immediately
follows from the first one. \qed

\begin{pro}\label{d-rel-b-fine} Let $\rd$ be a locally finite root system of type $B_I$ and $S_D\sub\rd$ be the subsystem of type $D_I$ introduced above.


(i) Suppose that $\dot\PP$ is a reflectable set for $S_D$ and there
exist $\da,\db\in\dot\PP$ with $\supp(\da)=\supp(\db)=\{i,j\}$ and
$(\da,\db)=0$. Then
$\dot\Pi:=(\dot\PP\setminus\{\da\})\cup\{\ep_i\}$ is a reflectable
set for $\rd.$

(ii) Suppose $\dot\PP\sub S_D$ and $i\in I$ are such that
$\dot\Pi:=\dot\PP\cup\{\ep_i\}$ is a reflectable set for $\rd.$ Let
$\da=r\ep_i+s\ep_j\in\dot\PP$ and set $\db:=r\ep_i-s\ep_j$. Then
$\dot\Pi':=\dot\PP\cup\{\db\}$ is a reflectable set for $S_D$.


\end{pro}

\proof
(i) By assumption, we have $\da=r\ep_i+s\ep_j$ and
$\db=r\ep_i-s\ep_j$ for some $i,j\in I$ and $r,s\in\{\pm1\}$.
 Now as $\db,\ep_i\in\dot\Pi,$ we have
$-\da=-r\ep_i-s\ep_j=w_{\ep_i}(\db)\in \w_{\dot\Pi}\dot\Pi$. This
gives  that $\dot\PP\sub \w_{\dot\Pi}\dot\Pi.$ Now as $\dot\PP$ is a
reflectable set for $S_D$, $S_D\sub \w_{\dot\Pi}\dot\Pi.$ Then for
any $k\in I,$ $\ep_k=w_{\ep_i-\ep_k}\ep_i\in \w_{\dot\Pi}\dot\Pi$.
Thus $\w_{\dot\Pi}\dot\Pi$ contains all nonzero roots of $\rd$ and
so $\dot\Pi$ is a reflectable set as required.

(ii) We first note that if $k\in I$ and $\dot\gamma\in
S_D^\times,$ then
\begin{equation}\label{nice}
\left\{
\begin{array}{ll}
w_{\ep_i}(\dgamma)=w_{\da}w_{\db}(\dgamma)&\hbox{if } i\in\hbox{supp}(\dgamma),\;\dgamma\neq\pm\da,\pm\db,\\
w_{\ep_i}(\dgamma)=\dgamma& \hbox{if }
i\not\in\hbox{supp}(\dgamma),\\
w_{\ep_i}(\dgamma)=\pm\db, & \hbox{if } \dgamma=\pm\da,\\
w_{\ep_i}(\dgamma)=\pm\da , & \hbox{if } \dgamma=\pm \db.
\end{array}\right.
\end{equation}
We must  show that $S_D^\times\sub \w_{\dot\Pi'}\dot\Pi'.$ Take
$\eta\in S_D^\times.$ Since $\dot\Pi$ is a reflectable set for
$\rd,$ there are $\da_1,\ldots,\da_t,\da'\in \dot\Pi$ such that
$w_{\da_1}\ldots w_{\da_t}(\da')=\eta$ which in turn implies that
$\da'\in\dot\PP.$ If for all $1\leq k\leq t,$ $\da_k\neq \ep_i,$
there is nothing to prove. Otherwise, suppose $1\leq k\leq t$ is such
that $\da_k=\ep_i$ and $\da_s\neq \ep_i$ for $k+1\leq s\leq t.$
Using (\ref{nice}), we can replace $w_{\da_k}\ldots w_{\da_t}(\da')$
with an element of $\w_{\dot\Pi'}\dot\Pi'.$ We do the same for other
reflections based on $\ep_i$ appearing in the expression and so we
get that  $\eta\in \w_{\dot\Pi'}\dot\Pi'.$\qed

\section{\bf Characterization of reflectable bases}\setcounter{equation}{0}\label{general-characterization}


In this section, we give a full characterization of reflectable sets
and reflectable bases for  tame irreducible reflection systems of
reduced types $X\neq E_{6,7,8}.$ We recall the definition of a coset
spanning set and a strong coset spanning set from Section
\ref{preliminaries}. Throughout this section, $(A,\fm,R)$ is  a tame
irreducible affine reflection system and $\w$ is its Weyl group. We
recall from Section \ref{preliminaries} that $\bar A=A/A^0$ and that
$\;^-:A\longrightarrow \bar A$ is the canonical epimorphism. We also
recall the map $\mathfrak{p}$ as in (\ref{p}).




\vspace{3mm} \noindent\underline{\bf Type $A_1$:} Considering
Theorem \ref{structure-thm}, we get that in this case
$R=(S+S)\cup(\rd^\times+S)$ in which $\rdot=\{\pm\da\}\cup\{0\}$ is
a finite  root system of type $A_1$ and $S$ is a pointed reflection
subspace of $A^0.$

\begin{thm}\label{finala1}
Suppose  $\Pi\sub\r^\times$  with $\la \Pi\ra=\la R\ra.$ Then
$\Pi$ is a reflectable set (resp. a reflectable base) for $R$ if
and only if $\Pi$ is  a strong coset spanning set (resp. a minimal
strong coset spanning set) for $R^\times$ in $2\la \r\ra$ with
respect to $\la\r\ra.$
\end{thm}

\proof Let $\Pi$ be a reflectable set. Then
$\la\Pi\ra=\la\w_\Pi\Pi\ra=\la R\ra$. Moreover, one can easily see
that
$$\w_\Pi\a\subseteq\a+2\la\Pi\ra\qquad(\a\in\Pi).$$
Therefore,
$$R^\times=\cup_{\a\in\Pi}\w_\Pi\a\subseteq\cup_{\a\in\Pi}[(\a+2\la \Pi\ra)\cap R^\times]=\cup_{\a\in\Pi}
[(\a+2\la R\ra)\cap R^\times]\subseteq R^\times.$$ This means that
$\Pi$ is a strong coset spanning set for $\r^\times$ in $2\la\r\ra$
with respect to $\la\r\ra.$

Next assume $\Pi$ is a strong coset spanning set for $\r^\times$ in
$2\la\r\ra$ with respect to $\la\r\ra.$ Then
\begin{equation}\label{temp11}
R^\times=\cup_{\a\in\Pi}[(\a+2\la R\ra)\cap R^\times].
\end{equation}
Since $\r^\times=\pm\da+S,$ using the  sign freeness of $\Pi$, we can (and do)
assume that  $\b-\mathfrak{p}(\b)=\da$ for all $\b\in\Pi.$ Since
$\da\in \rcross$ and $\la \Pi\ra=\la R\ra=\bbbz\da+\la S\ra,$ it
follows from (\ref{temp11}) that there exists $\d\in\la S\ra$ such
that $\a:=\da+2\d\in\Pi.$  Since $\Pi$ is a strong coset spanning
set for $\r^\times$ in $2\la\r\ra,$ it can be seen that
$\Pi_\a:=\Pi-\a$ satisfies
\begin{equation}\label{temp111}
\la \Pi_\a \ra=\la S\ra\andd S=\cup_{\sg\in\Pi_\a}(\sg+2\la
S\ra).\end{equation}

Next, we note that for $\zeta\in\Pi_\a,$ $-\a+\zeta=w_{\a}(\a+\zeta)$
which in turn implies that (by sign freeness)
\begin{equation}\label{sect2-2}
\begin{array}{c}
r\a+\zeta\in \w_{\Pi}\Pi\andd w_{\a+ r\zeta}\in\w_{\Pi};\\
r\in\{\pm1\},\;\zeta\in \Pi_\a.
\end{array}
\end{equation}

Also one can easily see that for $r,r_1,\ldots,r_n\in\{\pm1\}$ and
$\zeta,\zeta_1,\ldots,\zeta_n\in\Pi_\a,$ we have
\begin{equation*}
w_{\a+r_n\zeta_n}\ldots w_{\a+r_1\zeta_1}(r\a+\zeta)=(-1)^nr\a+\zeta-2rr_1\zeta_1+2rr_2\zeta_2+\cdots+2(-1)^nrr_n\zeta_n.
\end{equation*}
Now this together with (\ref{sect2-2}) and (\ref{temp111}) implies that
\begin{eqnarray*}
\w_{\Pi}\Pi&=&\pm\a+\Pi_\a+2\la\Pi_\a\ra\\
(\hbox{by} (\ref{temp111}))&=&
\pm\a+\Pi_\a+\la 2S\ra\\
&=&\pm\a+S\\
&=&R^\times
\end{eqnarray*}
which means that $\Pi$ is a reflectable set for $\r$.  This
 completes the proof.\qed

\vspace{3mm} \noindent\underline{\bf Type $B$:}


\begin{lem}\label{lemlem} Suppose that $\Pi\sub R^\times$ is a reflectable set for $\r,$ then $\Pi_{sh}$ is a
strong coset spanning set for $\r_{sh}$ in $\la\r_{lg}\ra$ with
respect to $\la\r\ra$ and $\Pi_{lg}$ is a coset spanning set for
$\r_{lg}$ in $2\la\r_{sh}\ra.$ Moreover  if $|I|=2,$ then
$\Pi_{lg}$ is a strong coset spanning set for $\r_{lg}$ in
$2\la\r_{sh}\ra$ with respect to $\la\r\ra.$ In other words,
$$\rsh=\bigcup_{\a\in\Pi_{sh}}[(\a+\la\rlg\ra)\cap\rsh],$$
$$\rlg=\bigcup_{\a\in\Pi_{lg}}[(\a+2\la \rsh\ra)\cap\rlg],\qquad(X=B_2)$$
and
$$\la\rlg\ra=\la\Pi_{lg}\ra+\la 2\rsh\ra.
$$

\end{lem}

\proof Using relations stated in Proposition \ref{structure-thm} and
the well known  facts about locally finite root systems, we have
\begin{equation}\label{newfine}
\la 2\rsh\ra+\la\rlg\ra\sub\la\rlg\ra\andd \la\rlg\ra+\la\rsh\ra\sub\la\rsh\ra.
\end{equation}
Now considering the different possibilities for $(\a,\b^\vee)$, $\a,\b\in R^\times$ and using (\ref{newfine}), one can check that
for $w\in\w$,
\begin{equation}\label{sect2-5}
w(\a)\in\a+\la R_{lg}\ra;\quad(\a\in\rsh,\;w\in\w)
\end{equation}
and
\begin{equation}\label{mmm1}
w(\b)\in\b+2\la R_{sh}\ra;\qquad (X=B_2,\;w\in\w,\;\b\in\rlg).
\end{equation}
Moreover, we note if $\a_1,\ldots,\a_m\in\r^\times$ and
 $\{\b_1,\ldots,\b_t\}=\{\a_1,\ldots,\a_m\}\cap R_{lg}$, then for $\a\in\r_{lg}$, we have
$$w_{\a_1}\cdots w_{\a_m}(\a)+2\la R_{sh}\ra=w_{\b_1}\cdots w_{\b_t}(\a)+2\la\r_{sh}\ra.
$$
Therefore, $\rlg=\w_\Pi\Pi_{lg}\sub\w_{\Pi_{lg}}\Pi_{lg}+2\la\rsh\ra$ and so
\begin{equation}\label{test100}
\la\rlg\ra\sub\la\Pi_{lg}\ra+2\la\rsh\ra\sub\la\rlg\ra.
\end{equation}

Now from (\ref{sect2-5}), (\ref{mmm1}), (\ref{test100}) and $\w\Pi=\r^\times$, we have
$$R_{sh}=\w\Pi_{sh}=\bigcup_{w\in\w,\;\a\in\Pi_{sh}}w(\a)\sub\bigcup_{\a\in\Pi_{sh}}[(\a+\la\rlg\ra)\cap\rsh]\sub\rsh,$$
$$ R_{lg}=\w\Pi_{lg}=\bigcup_{w\in\w,\;\a\in\Pi_{lg}}w(\a)\sub\bigcup_{\a\in\Pi_{lg}}[(\a+\la 2\rsh\ra)\cap\rlg]\sub\rlg\qquad(X=B_2),$$
and
$$ \la R_{lg}\ra=\la\Pi_{lg}\ra+2\la\rsh\ra.$$
Therefore the equalities in the statement hold. This completes the
proof.\qed

\vspace{3mm} \noindent\underline{\bf Type $B_2$:}
\begin{lem}\label{typeb2-1}
Consider a description
\begin{equation}\label{new-new}
R =(S+S)\cup(\rds+S)\cup(\rdl+L)
\end{equation}
for $\r$ as in Theorem \ref{structure-thm}. Let  $\RR_1:=\{\sg_i\mid
i\in J\}$ be a strong coset spanning set for $S$ in $\la L\ra$ and
$\RR_2:=\{\tau_j\mid j\in K\}$ be a strong coset spanning set for
$L$ in $\la 2S\ra$, such that
\begin{equation}\label{mal1}
\la\RR_1\ra+\la\RR_2\ra=\la S\ra.
 \end{equation}
 Assume that $0\in\RR_1$ and $0\in \RR_2$. For $i\in J$ and $j\in K$, pick $\da_i\in\rds$ and $\db_j\in\rdl$,
and set $\Pi_{sh}:=\{\da_i+\sg_i\mid i\in J\}$  and
$\Pi_{lg}:=\{\db_j+\tau_j\mid j\in K\}$. Then
$\Pi:=\Pi_{sh}\cup\Pi_{lg}$ is a reflectable set for $R.$
\end{lem}

\proof Since $\{\sg_i\mid i\in J\}$ is a strong coset spanning set
for $S$ in $\la L\ra$, we have by definition that  $S=\cup_{i\in
J}[(\sg_i+\la L\ra)\cap S]$ and $L=\cup_{j\in K}[(\tau_j+\la
2S\ra)\cap L]$. However, $S+L\sub S$ and $2S+L\sub L$, so using
(\ref{mal1}), we have
\begin{equation}\label{newtemp11}
\begin{array}{c}
S=\cup_{i\in J}(\sg_i+\la L\ra),\quad L=\cup_{j\in K}(\tau_j+\la 2S\ra)\andd 2\la\RR_1\ra+\la\RR_2\ra=\la L\ra.
\end{array}
\end{equation}

We must show that $\w_{\Pi}\Pi=\rcross.$  Since $0\in\RR_1$ and $0\in\RR_2$, we have $\rds\cap\Pi\not=\emptyset$ and $\rdl\cap\Pi\not=\emptyset.$
Then it is easy to see that $\rd^\times\sub\w_\Pi\Pi$ and so $\dot\w\sub\w_\Pi$. Therefore for $j\in J$ and $k\in K$,
$$\rds\pm\sg_j\sub\dot\w(\da_j\pm\sg_j)\sub\w_{\Pi}\Pi\andd\w_{\rds\pm\sg_j}\sub\w_{\Pi},
$$
and
$$\rdl\pm\tau_k\sub\dot\w(\db_k\pm\tau_k)\sub\w_{\Pi}\Pi\andd\w_{\rdl\pm\tau_k}\sub\w_{\Pi}.
$$
Next for $i,j\in J$ and $k,t\in K,$ we have
$$\rds\pm\sg_j\pm\tau_k\pm 2\sg_i\sub\dot\w\w_{\rds\pm\sg_i}\w_{\rdl\pm\tau_k}(\rds\pm\sg_j)\sub\w_\Pi\Pi,$$
and
$$\rdl\pm\tau_k\pm2\sg_j\pm 2\tau_t\sub\dot\w\w_{\rdl\pm\tau_t}\w_{\rds\pm\sg_j}(\rdl\pm\tau_k)\sub\w_\Pi\Pi.$$
For a fixed  $j\in J$, repeating this argument, we have
$$\rds\pm\sg_j+\la\tau_k\mid k\in K\ra+2\la\sg_i\mid i\in
J\ra\sub\w_\Pi\Pi,
$$
and
$$\rdl\pm\tau_k+2\la\tau_t\mid t\in K\ra+2\la\sg_i\mid i\in J\ra\sub\w_\Pi\Pi.
$$
Now using (\ref{newtemp11}) and (\ref{mal1}), we have
$$\rds+\sg_j+\la L\ra=\rds+\sg_j+\la\RR_2\ra+2\la\RR_1\ra\sub\w_\Pi\Pi,$$
and
$$\rdl+\tau_k+2\la S\ra=\rdl+\tau_k+2\la\RR_2\ra+2\la\RR_1\ra\sub\w_\Pi\Pi.
$$ Since $j\in J$ and $k\in K$ were chosen arbitrary, we get
from (\ref{newtemp11}) that
$$\rds+S\sub\w_\Pi\Pi\andd\rdl+L\sub\w_\Pi\Pi.$$
This completes the proof that $\Pi$ is a reflectable set.\qed

\begin{thm}\label{propro}
Suppose that $\Pi$ is  a subset of $R^\times$ with
$\la\Pi\ra=\la\r\ra.$ Then $\Pi$ is a reflectable set (resp.
reflectable base) for $\r$ if and only if $\Pi_{sh}=\Pi\cap R_{sh}$
is a strong coset spanning set (resp. a minimal strong coset
spanning set) for $\r_{sh}$ in $\la\r_{lg}\ra$ and $\Pi_{lg}=\Pi\cap
R_{lg}$ is a strong coset spanning set (resp. a minimal  strong
coset spanning set) for $\r_{lg}$ in $2\la\r_{sh}\ra,$ with respect
to $\la\r\ra.$
\end{thm}

\proof By Lemma \ref{lemlem}, it is enough to show  that if $\Pi$ is a subset of  $R^\times$ satisfying \begin{itemize}\item $\la\Pi\ra=\la\r\ra,$ \item $\Pi_{sh}$
is a strong coset spanning set for $\r_{sh}$ in $\la\r_{lg}\ra$ with respect to $\la\r\ra,$\item
$\Pi_{lg}$ is a strong coset spanning set for $\r_{lg}$ in
$2\la\r_{sh}\ra$ with respect to $\la\r\ra,$ \end{itemize} then it is a reflectable set.
So suppose  $\Pi\sub\r^\times$ satisfies the conditions mentioned
above. Therefore  $\Pi=\RR'_1\cup\RR'_2$, where $\RR'_1$ is a
strong coset spanning set for $\rsh$ in $\la\rlg\ra$ and $\RR'_2$
is a strong coset spanning set for  $\rlg$ in $2\la\rsh\ra$. Note
that $\bar{\Pi}$ contains at least one short root and one long
root and so is a reflectable set for $\bar{R}$ by Proposition
\ref{typeb2}. Therefore by Proposition \ref{end5}, it contains an
integral reflectable base $\PP$. Take a pre-image $\dot\Pi\sub\Pi$
of $\PP$ and construct a finite root system of type $B_2$, denoted again  by $\rd$,
 and $(S,L)$ as usual to get the description
(\ref{new-new}) for $R$ (see the proof of Theorem
\ref{structure-thm} for details). Since $\dot\Pi\sub\Pi\cap\rd$,
we have $\rd\cap\Pi_{sh}\not=\emptyset$ and
$\rd\cap\Pi_{lg}\not=\emptyset$. Let $\da\in \rd\cap\Pi_{sh},$ and
assume $\db\in\rd\cap\Pi_{lg},$  $\Pi_{sh}=\{\a_j\mid j\in J\}$
and $\Pi_{lg}=\{\b_t\mid t\in T\}$. We know that
$\mathfrak{p}(\a_j)\in S$ with $\mathfrak{p}(\da)=0$ and that
$\a_j-\mathfrak{p}(\a_j)-\da\in\la \rd_{lg}\ra$ for all $j\in J$.
Therefore for $j\in J,$ we have
$$\a_j+\la\rlg\ra=\da+\mathfrak{p}(\a_j)+\la\rlg\ra.$$
Now let $\tau\in S,$ then $\da+\tau\in \rsh=\cup_{j\in
J}[(\a_j+\la\rlg\ra)\cap\rsh].$ So
$\da+\tau\in\a_j+\la\rlg\ra=\da+\mathfrak{p}(\a_j)+\la\rlg\ra$, for
some $j\in J$. This gives $\tau\in\cup_{j\in
J}(\mathfrak{p}(\a_j)+\la L\ra).$ This means that
$\RR_1:=\{\mathfrak{p}(\a_j)\mid j\in J\}$ is a strong coset
spanning set  for $S$ in $\la L\ra$.  Using a similar argument, we
see that $\RR_2:=\{\mathfrak{p}(\b_t)\mid t\in T\}$ is
 a strong coset spanning set for $L$ in $2\la S\ra.$ We also
note that, as we have  already seen, $0\in\RR_1$ and $0\in \RR_2$.
Finally, since  $\la\Pi\ra=\la R\ra$, we have
$\la\RR_1\ra+\la\RR_2\ra=\la S\ra$. Therefore all conditions in the
statement of Lemma \ref{typeb2-1} hold and so $\Pi$ is a reflectable
set for $\r.$\qed

\vspace{5mm}

\noindent\underline{\bf Type $B_I$, $C_I$ $(|I|\geq 3)$:} We give
the proof for type $B$, the proof for type $C$ is analogous,
replacing the roles of short and long roots. So from now on, we
assume that we are in type $B$. We recall that in this case $L$ is a
lattice.

\begin{lem}\label{lem2-b>=3} Consider a description
\begin{equation}\label{new-new-new}
R =(S+S)\cup(\rds+S)\cup(\rdl+L)
\end{equation}
for $\r$ as in Theorem \ref{structure-thm}. Suppose that
$\RR:=\{\sg_j\mid j\in J\}$ is a strong coset spanning set for $S$
in $ L$ with $0\in\RR$, and $\MM:=\{\tau_k\mid k\in K\}$ is a
coset spanning set for $L$ in $\la 2S\ra.$ Suppose also that
$\dot\MM$ is a  coset spanning set  for $\rdl$ in $2\la\rds\ra$.
For each $j\in J$ and $k\in K$, pick $\da_j\in\rds$ and
$\db_k\in\rdl$ and set $\Pi:=\Pi_{sh}\cup\Pi_{lg}$ where
$\Pi_{sh}:=\{\da_j+\sg_j\mid j\in J\}$ and
$\Pi_{lg}:=\dot\MM\cup\{\db_k+\tau_k\mid k\in K\}$. Further,
suppose that
\begin{equation}\label{anny}
\la\RR\ra+\la\MM\ra=\la S\ra.
\end{equation}
Then $\Pi$ is a reflectable set for $R.$
\end{lem}

\proof First we note that since $\MM$ is a  coset spanning  set for
$L$ in $2\la S\ra,$ we have $\la \MM\ra+2\la S\ra=\la L\ra+2\la
S\ra=L$, therefore by (\ref{anny}),
\begin{equation}\label{anny2}
\la\MM\ra+2\la\RR\ra=L.
\end{equation}

Since $0\in\RR$, $\Pi$ contains at least a short root $\da\in\rd$.
Then we have from Proposition \ref{pro1-b>=3}(iii) that
$\{\da\}\cup\dot\MM$ is a reflectable set for the locally finite
root system $\rdot.$ So
\begin{equation}\label{4-b>=3}\rdot^\times\sub\w_\Pi\Pi\andd\dot\w\sub\w_\Pi.\end{equation}
From this, and the fact that $\a\in\w_\Pi\Pi$ if and only if $-\a\in\w_\Pi\Pi$,  we get
\begin{equation}\label{3-b>=3}
\rds\pm\sg_j=\dot\w(\da_j\pm\sg_j)\sub\w_\Pi\Pi;\qquad j\in J,
\end{equation}
and
\begin{equation}\label{5-b>=3}
\rdl\pm\tau_k=\dot\w(\db_k\pm\tau_k)\sub\w_\Pi\Pi;\qquad k\in K.
\end{equation}
Then for $i,j\in J$ and $k\in K$,
$$\rds\pm\sg_j\pm\tau_k\pm2\sg_i\sub\dot\w\w_{\rds\pm\sg_i}\w_{\rdl\pm\tau_k}(\rds\pm\sg_j)\sub\w_\Pi\Pi.$$
It follows, by repeating this argument, that for a fixed $j\in J$,
$$\rds+\sg_j+\la\MM\ra+2\la\RR\ra=\rds+\sg_j+\la\tau_k\mid k\in K\ra+2\la\sg_i\mid i\in J\ra\sub\w_\Pi\Pi.$$
Thus by (\ref{anny2}), $\rds+\sg_j+L\sub\w_\Pi\Pi$ for all $j\in J$. Now since $\cup_{j\in J}(\sg_j+L)=S$, we get
\begin{equation}\label{get1}
\rds+S\sub\w_\Pi\Pi.
\end{equation}

Next, we show that $\rdl+L\sub\w_\Pi\Pi$. Since $|I|\geq 3$, there
exist $\db_1,\db_2\in\rdl$ such that $(\db_1,\db_2^\vee)=\pm1$. This
together with (\ref{4-b>=3})-(\ref{5-b>=3}) shows that for $j\in J$
and $k,t\in K,$ we have
$$\rdl\pm\tau_k\pm2\sg_j\pm\tau_t\sub\dot\w\w_{\rdl\pm\tau_t}\w_{\rds\pm\sg_j}(\rdl\pm\tau_k)\sub\w_{\Pi}\Pi.
$$
By repeating this argument we get $$\rdl+\la\MM\ra+2\la\RR\ra=\rdl+\la\tau_k\mid k\in K\ra+2\la\sg_j\mid j\in J\ra\sub\w_\Pi\Pi,$$
and so by (\ref{anny2})
\begin{equation}\label{get2}
\rdl+L\sub\w_\Pi\Pi.
\end{equation}
Now (\ref{get1}) and (\ref{get2}) show that $\Pi$ is a reflectable set for $R$.\qed
\medskip

By a {\it $B_I$-reflectable data}, we mean a tuple
$$(\rdot,S,L,\dot\MM,\{\da_i\}_{i\in I},\{\db_k\}_{ k\in K},\RR,\MM)$$
in which the ingredients satisfy the conditions in the statement of
Lemma \ref{lem2-b>=3}.

\begin{thm}\label{nunny}
A subset $\Pi$ of $\r^\times$ with $\la\Pi\ra=\la\r\ra$ is a
reflectable set (resp. a reflectable base) for $\r$ if and only if
$\Pi_{sh}=\r_{sh}\cap\Pi$ is a strong coset spanning set (a
minimal strong coset spanning set) for $\r_{sh}$ in
$\la\r_{lg}\ra$ and $\Pi_{lg}=\r_{lg}\cap\Pi$ is a  coset spanning
set (a minimal    coset spanning  set)  for $\r_{lg}$ in
$2\la\r_{sh}\ra.$
\end{thm}

\proof By Lemma \ref{lemlem}, it is enough to prove that any
subset $\Pi$ of $\r^\times$ with $\la\Pi\ra=\la\r\ra$ such that
$\Pi_{sh}$ is a strong coset spanning set  for $\r_{sh}$ in
$\la\r_{lg}\ra$ and $\Pi_{lg}$ is a  coset spanning set for
$\r_{lg}$ in $2\la\r_{sh}\ra,$ is a reflectable set for $\r.$ So
take $\Pi$ to be such a subset of $\r^\times.$ Then
$\Pi=\RR_0\cup\MM_0,$ where $\RR_0$ is a strong coset spanning set
for $\r_{sh}$ in $\la\r_{lg}\ra$ and $\MM_0$ is a coset spanning
set for $\rlg$ in $2\la\rsh\ra.$ It follows that
$\bar{\Pi}={\bar\RR}_0\cup{\bar\MM}_0$ and that ${\bar\RR}_0$ is a
nonempty subset of $\bar\r_{sh}$ and ${\bar\MM}_0$ is a  coset
spanning set for ${\bar\r}_{lg}$ in $2\la{\bar\r}_{sh}\ra.$ Let
$\ab\in{\bar\RR}_0$. Then by Proposition \ref{pro1-b>=3},
$\{\ab\}\cup{\bar\MM}_0$ is a reflectable set for $\bar R$, so is
$\bar\Pi$. Thus by Corollary \ref{end5} and Proposition
\ref{end5}, it contains an integral reflectable base $\PP$. Take
$\dot\Pi\sub\Pi$ to be a pre-image of $\PP$ and construct $\rdot$
and $(S,L)$ as usual to get a description of $\r$ in the form
(\ref{new-new-new}), then $\dot\Pi$ is a reflectable base for
$\rdot$ and so by Proposition \ref{pro1-b>=3}, there is a short
root $\da$ and  a  coset basis $\dot\MM\sub\MM_0$  for $\rdl$ in
$\la \rds\ra$ such that $\dot\Pi=\{\da\}\cup\dot\MM.$


Next fix $\db\in\rdl$ and take $\d\in L,$ then
$\db+\d\in\la\r_{lg}\ra=\la\MM_0\ra+2\la\rsh\ra$ and so $\d\in
\mathfrak{p}(\la\MM_0\ra+2\la\rlg\ra)=\la
\mathfrak{p}(\MM_0\setminus\dot\MM)\ra+\la 2S\ra$. Thus $L=\la
\mathfrak{p}(\MM_0\setminus\dot\MM)\ra+2\la S\ra$,
which implies that $\mathfrak{p}(\MM_0\setminus\dot\MM)$ is a coset
spanning set for $L$ in $\la 2S\ra$.

Next, assume $\RR_0=\{\a_j\mid j\in J\}$,
$\MM_0\setminus\dot\MM=\{\b_k\mid k\in K\}$ and set
$$\RR:=\mathfrak{p}(\RR_0)=\{\sg_j:=\mathfrak{p}(\a_j)\mid j\in J\}\andd\MM:=\mathfrak{p}(\MM_0\setminus\dot\MM)=\{\tau_k:=\mathfrak{p}(\b_k)\mid k\in K\}.$$
Note that $0\in \RR\sub S$, as $\da\in\RR_0\cap\rd$. Now as in the
case of $B_2$, one sees that $\RR$  is a strong coset spanning set
for $S$ in $L$. It also follows that $\MM$ is a coset spanning set
for $L$ in $2\la S\ra$. Therefore, up to this point we have seen
that the tuple
$$(\rdot,S,L,\dot\MM,\{\da_i\}_{i\in I},\{\db_k\}_{ k\in K},\RR,\MM)$$
is a $B_I$ reflectable-data, provided that (\ref{anny}) holds. But
we know  that  $\la\Pi\ra=\la R\ra$, so
$\la\RR_0\cup\MM_0\ra=\la\rd\ra+\la S\ra$. Thus
$$\la S\ra=\mathfrak{p}(\la\RR_0\ra+\la\MM_0\ra)=\la \mathfrak{p}(\RR_0)\ra+\la \mathfrak{p}(\MM_0\setminus\dot\MM)\ra=
\la\RR\ra+\la\MM\ra,$$ so (\ref{anny}) holds. Now by Lemma
\ref{lem2-b>=3}, $\Pi$ is a reflectable set. The assertion in the
statement now follows immediately from Lemma \ref{lemlem}. \qed

\vspace{3mm}\noindent \underline{\bf Types $\dot A_I$$(|I|\geq 3)$
$D_I$$(|I\geq4|)$:} For the types under consideration, we have
\begin{equation}\label{fff1}R^\times=\rd^\times+S,
\end{equation}
where $S$ is a pointed reflection subspace of $A^0$ with $\la
S\ra=S$ and $\rd$ is a locally finite root system of the
corresponding type.

\begin{thm}\label{munny}
$\Pi\sub\r^\times$ is a reflectable set (resp. reflectable base) for
$R$ if and only if $\Pi$ is a  generating set (resp. minimal
generating set) for $\la\r\ra$.
\end{thm}

\proof We know that any reflectable set for $\r$ is a generating
set for $\la\r\ra,$ so we assume $\Pi\sub\r^\times$ is a
generating set for $\la\r\ra$ and show that it is a reflectable
set for $\r.$ By definition $\la\Pi\ra=\la R\ra$, and so
$\la\bar\Pi\ra=\la\bar R\ra$. Now by Proposition \ref{supp1},
$\bar\Pi$ is a reflectable set for $\bar{R}$, so by Propositions
\ref{a-new} and \ref{asim1-new}, it contains an integral
reflectable base $\PP.$ Let $\dot\Pi\sub\Pi$ be a pre-image of
$\PP$ and construct $\rd$ and $S$ in the usual way to get a
description of $R$ in the form (\ref{fff1}). Now $\dot\Pi$ is a
reflectable base for $\rd$ and so $\la\dot\Pi\ra=\la\rd\ra$. Set
$\MM=\Pi\setminus\dot\Pi$. Since $\la\rd\ra+S=\la R\ra=\la\Pi\ra$,
we have $S=\la \mathfrak{p}(\Pi)\ra=\la
\mathfrak{p}(\Pi\setminus\dot\Pi)\ra=\la \mathfrak{p}(\MM)\ra$.
Let $\MM=\{\a_j\mid j\in J\}$ and for $j\in J$ set
$\sg_j:=\mathfrak{p}(\a_j)$ and $\da_j=\a_j-\sg_j$. Since
$\dot\Pi\sub\Pi$, we have
$$\rd^\times=\w_{\dot\Pi}\dot\Pi\sub\w_{\Pi}\Pi\andd\dot\w\sub\w_\Pi.$$
Now for any $j\in J$, we have
$\rd^\times+\sg_j=\dot\w(\da_j+\sg_j)\sub\w_\Pi\Pi$. Since
$\a\in\w_\Pi\Pi$ if and only if $-\a\in\w_{\Pi}\Pi$, we have
$\rd^\times\pm\sg_j\sub\w_\Pi\Pi$. It is known that for the types
under consideration for any $\da\in\rd^\times$, there exists
$\db\in\rd^\times$ with $(\da,\db^\vee)=\pm1$. Therefore for any
$i,j\in J$, we have
$$\rd^\times\pm\sg_j\pm\sg_i\sub \dot\w\w_{\rd^\times\pm\sg_j}(\rd^\times\pm\sg_i)\sub\w_\Pi\Pi.$$
Repeating this argument, we get $\rd^\times+\la\sg_j\mid j\in J\ra\sub\w_\Pi\Pi$, so
$$R^\times=\rd^\times+S=\rd^\times+\la \mathfrak{p}(\MM)\ra=\rd^\times+\la\sg_j\mid j\in J\ra\sub\w_\Pi\Pi.$$
This shows that $\Pi$ is a reflectable set and  completes the
proof.\qed

 \vspace{3mm}\noindent \underline{\bf Types $F_4$,
$G_2$:}

\begin{thm}\label{g2f4}
Suppose that  $\PP$ is a subset of $\r^\times$ with
$\la\PP\ra=\la\r\ra.$ Then $\PP$ is a reflectable set (resp.
reflectable base) for $R$ if and only if $\PP_{sh}=\PP\cap\r_{sh}$
is a coset spanning set (resp. minimal coset spanning set) for
$\r_{sh}$ in $\la\r_{lg}\ra$ and $\PP_{lg}=\PP\cap\r_{lg}$ is a
coset spanning set  (resp. minimal coset spanning set) for $\r_{lg}$
in $\rho\la\r_{sh}\ra.$
\end{thm}
\proof Consider the description
\begin{equation}\label{f4g2-1}
R^\times=(\rds+S)\cup(\rdl+L)
\end{equation}
for $\r$ where $S$ and $L$ are as in Proposition
\ref{structure-thm} with $\rho=2$ for type $F_4$ and $\rho=3$ for
type $G_2$. We know that $\la S\ra=S$ and $\la L\ra=L$. let
$\PP\sub\r^\times$ be a reflectable set for $\r.$ Then by Lemma
\ref{general},
$$\rsh=\w_\PP\PP_{sh}=\bigcup_{\a\in\PP_{sh}}\w_\PP(\a)\sub\la\PP_{sh}\ra+\la\rlg\ra$$
and
$$\rlg=\w_{\PP}\PP_{lg}=\bigcup_{\a\in\PP_{lg}}\w_\PP(\a)\sub\la\PP_{lg}\ra+\la\rho\rsh\ra.$$
Thus we have   $\la\rsh\ra=\la\PP_{sh}\ra+\la\rlg\ra$ and
$\la\rlg\ra=\la\PP_{lg}\ra+\la\rho\rsh\ra.$ So by definition
$\PP_{sh}$ is a coset spanning set for $\rsh$ in $\la\rlg\ra$ and
$\PP_{lg}$ is a coset spanning set for $\rlg$ in $\la\rho\rsh\ra$.
Now as $\PP$ is a generating set for $\la\r\ra,$ we are done.

Next  Suppose that  $\Pi=\MM_1\cup\MM_2$, where $\MM_1$ is a coset
spanning set for $\r_{sh}$ in $\la\r_{lg}\ra$ and $\MM_2$ is a coset
spanning set for $\r_{lg}$ in $\la\rho\r_{sh}\ra$. Then
$\bar\Pi={\bar\MM}_1\cup{\bar\MM}_2$ such that $\bar\MM_1$ is a
coset spanning  set for $\rd_{sh}$ in $\la\rd_{lg}\ra$ and
$\bar\MM_2$ is a  coset spanning  set for $\rd_{lg}$ in
$\la\rho\rd_{sh}\ra$.  Thus by Propositions \ref{f4-locally} and
\ref{G2}, $\bar\Pi$ is a reflectable set for $\bar R$ and so by
Corollary \ref{end5}, it contains an integral reflectable base
$\PP$, for which we fix a pre-image $\dot\Pi\sub\Pi$, and we
construct $\rd$, $S$ and $L$ in the usual way to get a description
of $R$ in the form (\ref{f4g2-1}). Then $\dot\Pi$ is a reflectable
base for $\rd$. By Propositions \ref{f4-locally} and \ref{G2},
$\dot\Pi=\dot\MM_1\cup\dot\MM_2$ where $\dot\MM_1$ is a coset basis
for $\rds$ in $\la\rdl\ra$ and $\dot\MM_2$ is a coset basis for
$\rdl$ in $\la\rho\rds\ra)$. Let $\MM_1\setminus\dot\MM_1=\{\a_j\mid
j\in J\}$ and $\MM_2\setminus\dot\MM_2=\{\b_k\mid k\in K\}$. For
each $j\in J$, we have $\a_j=\da_j+\sg_j$ where
$\sg_j=\mathfrak{p}(\a_j)$ and $\da_j=\a_j-\mathfrak{p}(\a_j)$.
Similarly, for each $k\in K$, $\b_k=\db_k+\tau_k$ where
$\tau_k=\mathfrak{p}(\b_k)$ and $\db_k=\b_k-\mathfrak{p}(\b_k)$.
Take $\MM'_1=\mathfrak{p}(\MM_1\setminus\dot\MM_1)=\{\sg_j\mid j\in
J\}$ and $\MM'_1=\mathfrak{p}(\MM-2\setminus\dot\MM_2)=\{\tau_k\mid
k\in K\}$. One can check that $\MM'_1$ is a coset spanning set for
$S$ in $L$ and $\MM'_2$ is a coset spanning set for $L$ in  $\rho
S$. Since $\la\PP\ra=\la R\ra,$ it follows that
\begin{equation}\label{ggg1}
\la\MM'_1\ra+\la\MM_2'\ra=\la S\ra\andd
2\la\MM'_1\ra+\la\MM'_2\ra= L.
\end{equation}

Now as $\dot\Pi\sub\Pi$, we have
$$\rd^\times\sub\w_\Pi\Pi\andd\dot\w\sub\w_\Pi.$$
From this it follows that
$$\rds\pm\sg_j\sub\dot\w (\da_j\pm\sg_j)\sub\w_\Pi\Pi\andd\w_{\rds\pm\sg_j}\sub\w_{\Pi},\qquad(j\in J)$$
and
$$\rdl\pm \tau_k\sub\dot\w (\db_j\pm \tau_k)\sub\w_\Pi\Pi\andd\w_{\rdl\pm\tau_k}\sub\w_\Pi,\qquad(k\in K).$$
Next considering some basic facts on finite root systems of types
$F_4$ and $G_2$, we have for $i,j\in J$ and $t,k\in K$,
$$\rds\pm \sg_j\pm \sg_i\pm \tau_k\pm\tau_t\sub\dot\w\w_{\rds\pm \sg_i}\w_{\rdl\pm\tau_t}\w_{\rdl\pm\tau_k}(\rds\pm\sg_j)\sub\w_\Pi\Pi.
$$
Repeating this argument, and using (\ref{ggg1}), we obtain
$$\rds+S=\rds+\la\MM'_1\ra+\la\MM'_2\ra\sub\w_\Pi\Pi.$$
We also have for $i,j\in J$ and $t,k\in K$,
$$\rdl\pm\tau_k\pm\tau_t\pm\sg_i\pm\sg_j\sub\dot\w\w_{\rdl\pm\tau_t}\w_{\rds\pm\sg_i}\w_{\rds\pm\sg_j}(\rdl\pm\tau_k)\sub\w_\Pi\Pi.
$$
Again, repeating this argument and using (\ref{ggg1}), we get
$$\rdl+L=\rdl+2\la\MM'_1\ra+\la\MM'_2\ra\sub\w_\Pi\Pi.
$$
These all together show that $\Pi$ is a reflectable set. This
completes the proof.\qed

\end{document}